\documentstyle[11pt]{article}
\def\ns{\normalsize}                 \def\ds{\displaystyle}
\def\beq{\begin{equation}\label}     \def\eeq{\end{equation}}
\def\bea{\begin{eqnarray}\label}     \def\eea{\end{eqnarray}}
\def\ar{\begin{array}}               \def\er{\end{array}}
\def\la{\lambda}         \def\up{\stackrel}    \def\wt{\widetilde}
\def\rar{\rightarrow}      
 \def\wh{\widehat}
\renewcommand{\o}{\otimes}
\def\red{\hspace*{6mm}}
\def\al{\alpha}     \def\be{\beta}        \def\ga{\gamma}
            \def\B{{\cal B}}
\def\J{{\cal J}}    
\def\p{[p/\hbar]}     
     \def\dl{\delta}
\newcommand{\REF}[1]{$\hspace*{1pt}{}^{\rm #1}\,$}
\newtheorem{theorem}{Theorem}  
\newtheorem{prop}{Proposition} 
\newtheorem{defn}{Definition}  
\def\rem{{\it Remark:} \ }   \def\proof{{\it Proof}\ }
\begin{document}
\baselineskip 16pt
\thispagestyle{empty}

\vspace*{2.5mm}

August 1995 \hfill {q-alg/9508022}

\vspace*{1.6cm}  \begin{center}
{\Large \bf $\bf (T^*{\cal B})_q$,
 $\bf q$-analogue of model space  \\[2mm]
 and CGC generating matrices \ns $ ^{\dag}$ }
\end{center}

\vspace{1cm}

\begin{center}  {\bf A. G. Bytsko} $^{a)}$  \end{center}
\vspace*{0.1cm}

\begin{center} \begin{tabular}{c}
\it St.Petersburg Branch of Steklov Mathematical Institute \\
\it Fontanka 27, St.Petersburg 191011, Russia  \\
\end{tabular} \end{center}
 \vspace*{0.4cm}

\begin{center}  {\bf L. D. Faddeev} $^{b)}$  \end{center}
\vspace*{0.1cm}

\begin{center} \begin{tabular}{c}
\it St.Petersburg Branch of Steklov Mathematical Institute \\
\it Fontanka 27, St.Petersburg 191011, Russia  \\
    and                   \\
\it Research Institute for Theoretical Physics \\
\it P.O. Box 9 (Siltavuorenpenger 20C), University of Helsinki  \\
\it Helsinki SF-00014, Finland
\end{tabular} \end{center}

\vspace*{1.3cm}

\begin{center} \bf{ Abstract  } \end{center}
\red We study relations between the deformed cotangent bundle
$(T^*{\cal B})_q$ for the Borel subgroup $\cal B$ of a given simple
Lie group $G$, the quantum Lie algebra $\J_q$ associated with the
corresponding quantum group $G_q$ and the matrices generating
Clebsch-Gordan coefficients (CGC) for $\J_q$.
We reveal the connection of these objects to quantum analogue of
the model space $\cal M$ and $q$-tensor operators.

\vfill

\begin{flushleft}   \rule{6 cm}{0.05 mm} \\
$^{a)}$ {\footnotesize e-mail:$\;$ bytsko@pdmi.ras.ru } \\
$^{b)}$ {\footnotesize e-mail:$\;$ faddeev@pdmi.ras.ru } \\
$ ^{\dag}$ {\footnotesize
Published in Journal Math. Phys. {\bf 37}
%%, 6324-6348
(1996). } \\
\end{flushleft}

\newpage

\section{ \ns\bf INTRODUCTION. }

\red Among different representations of a given compact Lie group $G$
the model space ${\cal M}$ plays a distinguished role.
By definition,\REF{1} the model space is a direct sum of all
irreducible representations ${\cal H}_j$ with multiplicity one
\beq{0.0}
    {\cal M} = \sum\limits_j \oplus{\cal H}_j
\eeq
realized in some universal way. A most popular form of
${\cal M}$ is a space of holomorphic functions on  the Borel subgroup
$\cal B$ of complexified form of the group $G$. In this construction the
Borel  subgroup is considered as an affine space.

A study of model spaces provides a natural language for investigation of
physical models. For example, the popular model of $2$-dimensional quantum 
gravity, introduced by Polyakov,\REF{2} may be interpreted in terms of 
the model space of Virasoro algebra.\REF{3} 
A finite-dimensional quantum group with deformation parameter, depending
on the central charge, naturally appears in this context.

In the present paper, which was written with an intent to find new
applications of model space in modern mathematical physics, we discuss a
$q$-analogue of the model space related to $q$-deformed Lie group $G_q$. 
For this purpose we introduce and examine several "coordinatizations" of 
the quantum space $(T^*\B)_q$. As a by-product we obtain some generating 
matrices for the set of Clebsch-Gordan coefficients (CGC). To our 
knowledge this result is new even for the non-deformed case.

Throughout the paper we systematically and intentionally make use of the
$R$-matrix formalism, which we believe is the most convenient and powerful
tool to get explicit results in the domain of quantum groups.

To avoid the known difficulties with compact forms of quantum groups we 
adopt here a convention to work with complexified objects (groups, 
algebras) and their finite dimensional representations on a formal 
algebraic level. We also do not discuss subtleties arising in the 
case of $q$ being a root of unity.

Most of formulae given in this paper in $R$-matrix form have universal
validity. However, the concrete results are illustrated on the simplest 
example $G_q=SL_q(2)$. The generalization to other groups needs more
technical details such as an explicit structure of $R$-matrices and 
related objects.

Mentioned above "coordinatizations" of $(T^*\B)_q$ arise from two possible
decompositions of the matrix $L$ (in usual notations $L=L_+L_-^{-1}$, it
comprises all generators of the corresponding quantum Lie algebra) :
$$ L = U\,D\,U^{-1} \;\;\; {\rm and} \;\;\; L = A\,B\,A^{-1} \ , $$
where $D$ is a diagonal unimodular matrix, $U$ is a deformation of unitary
matrix, $A$ and $B$ are unimodular upper and lower triangular matrices. 
As we shall clarify below, the matrices $A$ and $B$ admit a natural
interpretation as the coordinates in the base and in the fiber of
$(T^*\B)_q$, whereas entries of the matrix $U$ will be shown to
provide basic shifts on
the model space $\cal M$ and generate $q$-analogues of Clebsch-Gordan
coefficients for the quantum group $G_q$. The explicit connection between
$U$ and $(A,B)$ will be demonstrated on the example of  $SL_q(2)$.

It should be mentioned that an object like the matrix $U$ appeared first
in Refs.4,5 (later it was used also in Ref.6), where it was interpreted
as a "chiral" component of the quantum group-like element $g$.
In the present paper we give another interpretation and application of
the matrix $U$ in the context of a model space.

Let us briefly describe the contents of the present paper. In the Sec.$\!$
II the definition of the cotangent bundle for a quantum group is reminded.
Next we introduce an object of especial interest for us -- the algebra
$\cal U$ generated by the entries of the matrix $U$ which diagonalizes the
coordinate in a fiber of $(T^*G)_q$. We derive explicit relations for this
algebra in the case of $G=SL(2)$.

In the Sec.$\!$ III we consider a non-deformed limit ($q=1$) of the algebra
$\cal U$ and construct an explicit representation.
For the case of $G=SL(2)$ we show that the matrix
$U_0$  generates Clebsch-Gordan coefficients (CGC) for the corresponding
non-deformed Lie algebra. The Borel subgroup $\cal B$ and the space $T^*\B$
naturally appear here.  Finally, we discuss a connection of our results
with the  Wigner-Eckart theorem.

In the Sec.$\!$ IV we construct representations of the algebra
$\cal U$ (for $q\neq 1$) for the case of $SL(2)$ in two different ways.
The first one uses the language of $q$-oscillators.
The second is based on explicit realization
of $(T^*\B)_q$ and hence involves a notion of quantum model space.
Here we show that the matrix $U$ is a "generating matrix" for CGC for
deformed Lie algebra.  We also give some comments on the generalized
version of the Wigner-Eckart theorem.

\section{\ns\bf $\bf (T^*G)_q$ AND RELATED OBJECTS.}
\setcounter{equation}{0}

\red There exist three symplectic manifolds (from the physical point
of view they are phase spaces) naturally related to a given Lie group
$G$ and its Lie algebra ${\cal J}$:

$\langle 1 \rangle$\ \  $T^*G$ -- the cotangent bundle for the group $G$;

$\langle 2 \rangle$\ \  $T^*\B$ --
the cotangent bundle for the Borel subgroup $\B$;

$\langle 3 \rangle$\ \  $\cal O$ --
an orbit of the co-adjoint action of $G$ on ${\cal J}^*$.

For instance, in the case of $G=SL(2)$ (which will be our main example)
these spaces are six-, four-, and two-dimensional, correspondingly.

The method of geometric quantization\REF{7} provides a 
representation theory for $\langle 1 \rangle$, $\langle 2 \rangle$ and
$\langle 3 \rangle$. Turning from classical to quantum 
groups, one can try to construct a representation theory for the deformed 
analogues of these manifolds. In the present paper we shall deal with
deformations of the spaces $\langle 1 \rangle$ and $\langle 2 \rangle$.

\subsection*{\red \ns\bf A. Description of $\bf (T^*G)_q$. }
\red Let $G_q$ be a deformation of the Lie group $G$ and ${\cal J}_q$ be 
a deformation of the corresponding Lie algebra ${\cal J}$.
The deformed cotangent bundle $(T^*G)_q$ is a non-commutative manifold,
i.e., according to the ideology developed by A.Connes,\REF{8} its 
coordinates are (non-commuting) generators of some associative algebra. 
A point on this manifold is parameterized by the pair $(g,L)$, where
$g\in G_q$ is a coordinate in the base of the bundle, and $L$ is a
coordinate in a fiber.

The structure of $(T^*G)_q$ is defined via commutation relations between the
coordinates in the base and in a fiber. An appropriate $R$-matrix form of
these relations was proposed in Ref.5:
\beq{1.1}
 R_{\pm} \up{1}{g} \,\up{2}{g} \,=\, \up{2}{g} \, \up{1}{g} R_{\pm} \,,
\eeq
\beq{1.2}
R_- \up{1}{g} \up{2}{L} \,=\, \up{2}{L} R_+
\up{1}{g} \,,
\eeq
\beq{1.3}
 \up{1}{L} R_-^{-1} \up{2}{L} R_- \,=\, R_+^{-1}
\up{2}{L} R_+ \up{1}{L} \,.
\eeq
Here and below we use the formalism developed in Ref.9, i.e., objects 
like $g$ and $L$ are considered as matrices (say, $L\in \J_q\otimes V$, 
where $V$ stands for auxiliary space). We use the standard notations 
for tensor products: 
$\up{1}{L}\up{}{=L\otimes I}\in\J_q\otimes V\otimes V$,  etc.

Let us take the parameter $q$, which appears in the theory of 
quantum groups, in the following form
\beq{0}
q=e^{\gamma\hbar} ,
\eeq
where $\hbar$ is the Planck constant (the parameter of quantization)
and $\ga$
is the deformation parameter. In physical applications it is most natural
to suppose that $\ga$ is either pure real ($q$ belongs to the real axis)
or pure imaginary ($q$ belongs to the unit circle at the complex plane).

The second form of $q$ is typical for the WZW theory.\REF{4,6,10}
{}For $|q|=1$ we suppose also that $q$ is not a root of unity.
It should be mentioned that for both variants of choice of $\ga$
in (\ref{0}) the definition of $q$-number
\beq{0.2}
[x] \equiv \frac{q^x - q^{-x}}{q-q^{-1}}
\eeq
is invariant with respect to complex conjugation of $q$, i.e.
$\overline{[x]}=[\overline{x}]$. This property
becomes important if one discusses involutions of deformed Lie algebras.

\begin{defn}\label{DL}
The algebra $\cal L$ is an associative algebra generated by entries of 
the  matrix $L$ which obeys relation (\ref{1.3}).
\end{defn}

An important fact --  the connection of algebra $\cal L$ with the 
corresponding quantum Lie algebra $\J_q$ was established in Ref.11
in the following form.
\begin{prop}\label{P4}
Let matrices $L_+$, $L_-$ obey the following exchange relations
\beq{L2}
R_{\pm} \up{1}{L_+} \up{2}{L_+} = \up{2}{L_+} \up{1}{L_+} R_{\pm} ,\;\;
R_{\pm} \up{1}{L_-} \up{2}{L_-} = \up{2}{L_-} \up{1}{L_-} R_{\pm}, \;\;
R_{+} \up{1}{L_+} \up{2}{L_-} = \up{2}{L_-} \up{1}{L_+} R_{+} \;\; .
\eeq
Then the matrix $ L= L_+ \, L_{-}^{-1}$ satisfies the relation (\ref{1.3}).
\end{prop}
This statement implies that the algebra $\cal L$ is isomorphic 
(up to some technical details which we do not discuss here) to 
corresponding quantum Lie algebra $U_q(\J)$ [which is defined by (\ref{L2}), 
see, e.g., Ref.9].

Consider now the relations (\ref{1.1})-(\ref{1.3})
for $g$ and $L$ being $2\times 2$ matrices
\beq{g,L}
g=\left(\ar{cc} g_1&g_2\\g_3&g_4\er\right),\;\;\;\; 
L=\left(\ar{cc}A&B\\C&D\er\right)
\eeq
and the $R$-matrices taken in the form
\beq{1.4}
R_{+}=q^{-1/2}\left(\ar{cccc}
q & 0 & 0 & 0 \\0 & 1 & \omega & 0 \\0 & 0 & 1 & 0 \\0 & 0 & 0 & q \er
\right) ,\; \omega \equiv q-q^{-1}; \;\;\;\;\;\;   R_{-}=PR_{+}^{-1}P
\eeq
($P$ denotes the permutation operator: $P\!\up{1}{g}\!P\!=\,\up{2}{g}$,
etc.). In this case (\ref{1.1})-(\ref{1.3}) define the cotangent bundle for 
the quantum group $G_q=GL_q(2)$; each of $R$-matrix equations (\ref{1.1}) 
and  (\ref{1.3}) is equivalent to six independent relations:
\beq{1.35}
\ar{c} q\,g_1\,g_2= g_2\,g_1,\;\; q\,g_1\,g_3= g_3\,g_1,\;\; 
q\,g_2\,g_4= g_4\,g_2,\;\;  q\,g_3\,g_4= g_4\,g_3,\;\; \\ \\ 
g_2\,g_3= g_3\,g_2,\;\;\;\; g_1\,g_4 - q^{-1}\, g_4\,g_1 = 
-\omega \,g_2\,g_3 ; \er
\eeq
and
\beq{1.5}
\ar{c} [A,B]=-q^{-1}\omega BD,\;\;\; [A,C]=q^{-1}\omega DC,\;\;\;[A,D]=0, \\
\\ CD=q^2DC,\;\;\;BD=q^{-2}DB,\;\;\;[B,C]=q^{-1}\omega D (D-A). \er
\eeq
The equation (\ref{1.2}) gives the following relations
\beq{gL}
\ar{ll} g_1 A = q A g_1 + \omega B g_3,\; & g_1 B = B g_1, \\ [1mm]
g_2 A = q A g_2 + \omega B g_4,\; & g_2 B = B g_2, \\ [1mm]
g_3 A = q^{-1} A g_3 + \omega g_1 C,\; & g_3 B = B g_3
+\omega g_1 D, \\ [1mm]
g_4 A = q^{-1} A g_4 + \omega g_2 C,\; & g_4 B = B g_4
+\omega g_2 D, \\ [1mm]
g_1 C = C g_1 + q^{-1}\omega D g_3,\; & g_1 D = q^{-1}\, D g_1, \\ [1mm]
g_2 C = C g_2 + q^{-1}\omega D g_4,\; & g_2 D = q^{-1}\, D g_2, \\ [1mm]
g_3 C = C g_3, & g_3 D = q D g_3, \\ [1mm]
g_4 C = C g_4, & g_4 D = q D g_4. \er
\eeq

Next, let us recall the well-known statement (see, e.g., Ref.9):
\begin{prop}\label{P1}
The algebra generated by the entries of the matrix $g$ obeying (\ref{1.35})
possesses the central element (``deformed determinant'')
\beq{1.65}
{\det}_q g = g_1\,g_4-q^{-1}g_2\,g_3 .
\eeq
\end{prop}

Similarly, for the algebra $\cal L$ in the case of $GL_q(2)$ one can
check the following.
\begin{prop}\label{P2}
The algebra with generators $A$, $B$, $C$, $D$ obeying (\ref{1.5}) possesses 
two central elements:
\beq{1.6}
K_1 = qA + q^{-1} D, \;\;\; K_2 = q^{-1} AD - q BC.
\eeq
\end{prop}

Finally, using the commutation relations (\ref{gL}), one can check that
\begin{prop}\label{P3}
The operators ${\det}_q g$ and $K_2$ commute with all entries of the 
matrices $g$ and $L$. 
\end{prop}
This implies that, fixing values of ${\det}_q g$ and $K_2$, one gets a
certain subalgebra of the algebra defined by (\ref{1.35})-(\ref{gL}).
\begin{defn}\label{D1}
Relations (\ref{1.35})-(\ref{gL}) for 
${\det}_q g = 1$ and $ K_2 = {\rm const}$ 
define the cotangent bundle for the quantum group $G_q=SL_q(2)$.
\end{defn}

Let us underline that the above definitions and statements can be easily
generalized, say to the case of $SL_q(N)$.

In our case the algebra $\cal L$ is isomorphic to the quantum Lie algebra
${\cal J}_q=U_q(sl(2))$ (introduced first in Ref.12) which is defined by
the relations
\beq{L3}
[l_+ , l_- ] = \frac{q^{2l_3}-q^{-2l_3}}{q-q^{-1}}\equiv [2l_3],\;\;\;\;
q^{l_3} \, l_{\pm}  = q^{\pm 1} l_{\pm} \, q^{l_3}
\eeq
and the matrices $L_{\pm}$ can be chosen as follows:
\beq{L4}
L_+ = \left( \ar{cc} q^{l_3} & \omega q^{1/2} l_- \\ 0 &
q^{-l_3} \er \right),\;\;\;\;
L_- = \left( \ar{cc} q^{-l_3} & 0 \\
-\omega q^{-1/2} l_+ & q^{l_3} \er \right) .
\eeq

Note that the matrix $L$ in the Proposition \ref{P4} is defined only up to
a scaling factor. Thus, for $L_+$, $L_-$ given in (\ref{L4}), we may choose 
$L$ as follows
\beq{L5}
L= q^2 L_+ L_{-}^{-1} =
\left(\ar{cc} q C - q^{-2 l_3} &  q^{5/2} \omega l_-
q^{-l_3} \\ q^{-1/2}\omega l_+ q^{-l_3} & q^2 q^{-2l_3} \er \right).
\eeq
Here $C$ stands for the Casimir operator of $U_q(sl(2))$ :
\beq{L8}
 C = {\omega}^2 l_- l_+ + q^{2l_3+1} + q^{-(2l_3+1)} = q^{2\wh{j}+1} +
 q^{-(2\wh{j}+1)} ,
\eeq
where $\wh{j}$ is the operator of spin.

According to Proposition \ref{P4}, the matrix (\ref{L5}) satisfies 
(\ref{1.3}). Therefore, it provides a (fundamental) representation of the 
algebra ${\cal L}$ for $U_q(sl(2))$. In this representation the central 
elements (\ref{1.6}) are given by
\beq{L6}
K_1 = q^2 C, \;\;\; K_2 = q^3,
\eeq
Note that the scaling factor $q^{2}$ introduced in (\ref{L5}) has changed
the values of $K_1$ and $K_2$. The choice of such normalization in (\ref{L5})
will be explained later.

\subsection*{\red \ns\bf B. Connection with quantum $\bf 6j$-symbols. }
\red  Let us remind the theorem which describes an important property of 
the algebra ${\cal L}$ for $U_q(sl(2))$ (this statement
appeared first in Ref.5).
\begin{theorem}
 Let $D\equiv D(p)$ be the unimodular diagonal matrix 
\beq{1.75}
 D = \left( \ar{cc}
q^{p/\hbar} &  \\ & q^{-p/\hbar} \er \right) ,
\eeq
and let 2$\times$2 matrix $U$ satisfy the following exchange relations 
\beq{1.8}
\up{1}{D}\, \up{2}{U} = \up{2}{U} \up{1}{D} \sigma,\;\;
\; \up{2}{D}\, \up{1}{U} = \up{1}{U} \up{2}{D}
\sigma, \;\;\; \sigma= {\rm diag} (q^{-1},q,q,q^{-1}),
\eeq
\beq{1.9}
R_+ \up{1}{U} \, \up{2}{U} = \up{2}{U} \,\up{1}{U} {\cal R}_+ (p), \;\; 
R_- \up{1}{U}\, \up{2}{U} = \up{2}{U}\, \up{1}{U} {\cal R}_- (p),
\eeq
where $R_{\pm}$ are the standard $R$-matrices (\ref{1.4}) and
\beq{1.10}
{\cal R}_+ (p)= P{\cal R}_-^{-1} (p) P = q^{-1/2} \left( \ar{cccc}
q &  &  &   \\
 & \ds \frac{\sqrt{ [p/\hbar +1] [p/\hbar -1] }}{[p/\hbar]} &     \ds
\frac{q^{p/\hbar}}{[p/\hbar]} &  \\  &  &  \\
 & \ds -\frac{q^{-p/\hbar}}{[p/\hbar]} & \ds
\frac{\sqrt{ [p/\hbar +1] [p/\hbar -1] }}{[p/\hbar]} &  \\
 &  &  & q  \er \right),
\eeq
(here $[x]$ denotes a "q-number" (\ref{0.2})).
Then matrix $L$ constructed by means of the similarity transformation
\beq{1.7}
L = U D\, U^{-1}, \;\;\;
\eeq
satisfies the relation (\ref{1.3}) and therefore its entries generate an
algebra ${\cal L}$ isomorphic to $U_q(sl(2))$.
\end{theorem}

The proof is given in appendix A. It makes use of the identity
\beq{RD}
{\cal R}_- (p) = (\up{1}{D})^{-1} {\cal R}_+ (p) \sigma \up{1}{D} .
\eeq

\rem A consequence of (\ref{1.8}) is the commutativity of $L$ and $D$
\beq{1.12}
\up{1}{L} \, \up{2}{D} = \up{2}{D} \, \up{1}{L},
\eeq 
which implies that $p$ commutes with all elements of $\cal L$. Later 
we shall interpret $p$ as the operator of spin.

\rem Properly  generalizing objects which enter the Theorem 1, one 
can extend this theorem to the case of any quantum semisimple Lie 
algebra.\REF{13} 
In particular, the matrix $D$ for $U_q(sl(N))$ is found to be :
$D(\vec{p}\,)\,=\, const\cdot q^{\vec{H}\,\o\,\ds\vec{p} }$,\ 
where $\vec{p}$ consists of the operators corresponding to components 
of the weight vector (i.e., on each irreducible representation they are
multiples of unity) 
and $\vec{H}$ consists of the generators $H_i$ of the Cartan subalgebra.
An explicit form of ${\cal R}(p)$ for $U_q(sl(N))$ was obtained in Ref.14. 

\rem The matrix ${\cal R}(p)$ obeys the deformed Yang-Baxter 
equation,\REF{14-16,5} which can be written, for example, as follows :
\beq{1.13}
\up{1}{Q} \up{23}{\cal R}_+ (p) (\up{1}{Q})^{-1} \up{13}{\cal R}_+ (p)
\up{3}{Q} \up{12}{\cal R}_+ (p) (\up{3}{Q})^{-1} =
\up{12}{\cal R}_+ (p) \up{2}{Q} \up{13}{\cal R}_+ (p) (\up{2}{Q})^{-1}
\up{23}{\cal R}_+ (p),
\eeq
where for $\J_q=U_q(sl(2))$ the matrix 
$Q = \biggl(\! \ar{cc} e^{i\xi} & \\ &\! e^{-i\xi} \er\!\biggr)$
contains an extra variable $\xi$, conjugated with $p$:
\beq{1.15}
[p,\xi] = -i \hbar, \;\;\;\;\; q^{p/\hbar}\,e^{i\xi} =
q\,e^{i\xi}\,q^{p/\hbar}.
\eeq
This variable $\xi$ belongs to the algebra $\cal U$ but does not enter matrix
$L$. An explicit expression for $\xi$ will be given below. The general form
of $Q$ for $U_q(sl(N))$ can be easily found\REF{13} : 
$ Q\,=\,e^{ i\,\vec{H}\,\o\,\ds\vec{\xi} }$,  where components of $\xi_i$ 
are operators conjugated to $p_i$ : $[p_j,\xi_k]=-i\hbar\,\dl_{jk}$.

The matrix ${\cal R}(p)$ was discussed in physical literature in different 
contexts.  In particular, it plays significant role in studies of quantum 
Liouville\REF{15,16} and WZW\REF{4-6} models; its relation to 
Calogero-Moser model was recently discussed in Ref.17. 
But for us more important fact is a connection of ${\cal R}(p)$ with the 
quantum $6j$-symbols:
the entries of (\ref{1.10}) calculated on irreducible representations 
coincide (up to some normalization) with the values of some $6j$-symbols 
for $U_q(sl(2))$ (exact formulae are given in Ref.18, generalizations 
are discussed in Ref.13). This connection 
allows to assume that objects like the matrix $U$ should be interpreted in 
terms of Clebsch-Gordan coefficients (CGC). Below we demonstrate that $U$
is indeed a ``generating matrix'' for CGC and clarify
its  relation to $(T^*\B)_q$.

\subsection*{\red \ns\bf C. Algebra ${\cal U}$. }
\begin{defn}\label{DU}
 The algebra $\cal U$ is an associative algebra generated
by entries of matrix 
$U = \left( \ar{cc} U_1 & U_2 \\ U_3 & U_4 \er\right)$ and the operator
$p$ such that relations (\ref{1.75})-(\ref{1.10}) hold. 
\end{defn}
\rem For simplicity we restricted our consideration to the case of 
$\cal U$ associated with $U_q(sl(2))$. Let us stress that the case
of  $\cal U$ associated with $U_q(sl(N))$ can be studied
similarly but it will involve more technical details.
On the other hand, it might be rather cumbrous to obtain
exact formulae for $\cal U$ associated with $U_q(\J)$ in the case 
of $\J$ being generic semisimple Lie algebra.

Let us give an explicit form of the defining relations (\ref{1.9}) :
\beq{2.1}
U_1 U_3 = q^{-1} U_3 U_1, \;\;\;\;\; U_2 U_4 = q^{-1} U_4 U_2,
\eeq
\beq{2.2}
U_1 U_2 = U_2 U_1 \sqrt{\frac{[p/\hbar -1]}{[p/\hbar +1]}}, \;\;\;
U_3 U_4 = U_4 U_3 \sqrt{\frac{[p/\hbar -1]}{[p/\hbar +1]}}
\eeq
\beq{2.3}
U_1 U_4 = U_4 U_1 \frac{\sqrt{[p/\hbar +1] [p/\hbar -1]}}
{[p/\hbar]} - U_3 U_2 \frac{q^{p/\hbar}}{[p/\hbar ]}
\eeq
\beq{2.4}
U_3 U_2 = U_2 U_3 \frac{\sqrt{[p/\hbar +1] [p/\hbar -1]}}{[p/\hbar]}
- U_1 U_4 \frac{q^{-p/\hbar}}{[p/\hbar ]}
\eeq
The rest of the relations contained in (\ref{1.9}) are not independent and
can be deduced from (\ref{2.1})-(\ref{2.4}).

Additionally, from (\ref{1.8}) one gets
\beq{2.5}
\ar{ll} q^{p/\hbar} \, U_1 = q^{-1} \, U_1 \, q^{p/\hbar}, \; & q^{p/\hbar} \,
U_2 = q \, U_2 \, q^{p/\hbar},\\ & \\ q^{p/\hbar} \, U_3 = q^{-1} \, U_3 \, 
q^{p/\hbar},\;& q^{p/\hbar} \, U_4 = q \, U_4 \, q^{p/\hbar}.\er
\eeq

Thus, relations (\ref{2.1})-(\ref{2.5}) describe the algebra ${\cal U}$. 
Using them, one may verify the following statement: 
\begin{prop}
A central element of ${\cal U}$ is given by the "deformed" determinant of
the matrix $U$:
\beq{2.6}
{\rm Det} U \equiv U_1 U_4 \sqrt{\frac{[p/\hbar +1]}{[p/\hbar ]}}
- U_2 U_3 \sqrt{\frac{[p/\hbar -1]}{[p/\hbar]}}=
q U_4 U_1 \sqrt{\frac{[p/\hbar -1]}{[p/\hbar]}} -
q U_3 U_2 \sqrt{\frac{[p/\hbar +1]}{[p/\hbar ]}} .
\eeq
\end{prop}

For fixed value of ${\rm Det} U$ the algebra ${\cal U}$ contains only four 
independent generators. In classical limit $(\hbar = 0)$ they become the 
coordinates on 4-dimensional phase space.

For further discussion it is convenient to introduce new variables instead
of $U_i$:
\beq{2.7}
\wh{U}_i = U_i \sqrt{[p/\hbar]}.
\eeq

The coordinates $\{p$, $\wh{U}_i\}$ form a new set of generators of the algebra
${\cal U}$. The commutation relations (\ref{2.1})-(\ref{2.5}) rewritten
in terms of the new generators acquire a simpler form:
\beq{2.8}
\wh{U}_1 \wh{U}_3 = q^{-1} \wh{U}_3 \wh{U}_1, \;\;
\wh{U}_2 \wh{U}_4 = q^{-1} \wh{U}_4 \wh{U}_2, \;\;
\wh{U}_1 \wh{U}_2 = \wh{U}_2 \wh{U}_1, \;\;
\wh{U}_3 \wh{U}_4 = \wh{U}_4 \wh{U}_3
\eeq
\beq{2.10}
\wh{U}_1 \wh{U}_4 = \wh{U}_4 \wh{U}_1 \frac{[p/\hbar +1]}{[p/\hbar]} -
\wh{U}_3 \wh{U}_2 \frac{q^{p/\hbar}}{[p/\hbar ]},
\eeq
\beq{2.11}
\wh{U}_3 \wh{U}_2 = \wh{U}_2 \wh{U}_3 \frac{[p/\hbar +1] }{[p/\hbar]} -
\wh{U}_1 \wh{U}_4 \frac{q^{-p/\hbar}}{[p/\hbar ]},
\eeq
\vspace*{1mm}
\beq{2.9}
\ar{ll} q^{p/\hbar} \, {\wh U}_1 = q^{-1} \, {\wh U}_1 \, q^{p/\hbar}, \;&
q^{p/\hbar} \, {\wh U}_2 = q \, {\wh U}_2 \, q^{p/\hbar}, \\ & \\
q^{p/\hbar} \, {\wh U}_3 = q^{-1} \, {\wh U}_3 \, q^{p/\hbar},\;&
q^{p/\hbar} \, {\wh U}_4 = q \, {\wh U}_4 \, q^{p/\hbar}. \er
\eeq

The central element (\ref{2.6}) in new variables looks as follows
\beq{2.12}
{\rm Det} U \equiv ( \wh{U}_1 \wh{U}_4 - \wh{U}_2 \wh{U}_3)\frac{1}
{[p/\hbar ]} = ( \wh{U}_4 \wh{U}_1 - \wh{U}_3 \wh{U}_2 )\frac{q}{[p/\hbar]}.
\eeq

The explicit form of the matrix inverse to $\wh{U}$, which we shall need
later, is
\beq{2.125}
 \wh{U}^{-1} = \frac{1}{{\rm Det}U} \left(\ar{cc} \wh{U}_4 & -q \wh{U}_2 \\
 - \wh{U}_3 & q \wh{U}_1 \er\right)\frac{1}{[p/\hbar]}.
\eeq

Finally, from (\ref{2.7}) we conclude that the expression (\ref{1.7}) for
the matrix $L$ looks similarly in terms of new matrix $\wh{U}$:
\beq{2.130}
L \,=\, U\,D\,U^{-1} \,=\, \wh{U} \, D \, \wh{U}^{-1} \ .
\eeq

\section{ \ns\bf NON-DEFORMED CASE. }
\setcounter{equation}{0}
\subsection*{\red \ns\bf A. Representation of algebra ${\cal U}_0$. }
\red First, we consider the limit $\gamma\rightarrow 0$, $\hbar \neq 0$
(note that $q$-numbers turn into ordinary numbers),i.e., here we deal
with a well understood situation -- the representation theory of $SL(2)$.
An investigation of this simple non-deformed case will make further
results more transparent.
 
Let us denote the corresponding limit algebra as ${\cal U}_0$. The defining 
$R$-matrix relations (\ref{1.9}) now degenerate to
\beq{2.12b}
\up{1}{U_0}\, \up{2}{U_0} = \up{2}{U_0} \, \up{1}{U_0} {\cal R}_{\pm}^0 (p),
\eeq
where
\beq{1.12c}
{\cal R}_{+}^0 (p) = {\cal R}_{-}^0 (p) = \left( \ar{cccc} 1 &  &  & \\
& \ds\frac{\sqrt{ (p/\hbar +1) (p/\hbar -1) }}{(p/\hbar)} & 
\ds\frac{\hbar}{p} & \\ &  &  \\ & -\ds\frac{\hbar}{p} & 
\ds\frac{\sqrt{(p/\hbar +1) (p/\hbar -1)}}{(p/\hbar)}
& \\ &  &  & 1 \er\right).
\eeq

The analogues of relations (\ref{2.8})-(\ref{2.9}) for ${\cal U}_0$ are (from
now on we omit the index $0$ for the generators of ${\cal U}_0$)
\beq{2.13a}
 p \wh{U}_1 = \wh{U}_1 (p-\hbar), \;\; p \wh{U}_2 = \wh{U}_2 (p+\hbar),
\;\; p \wh{U}_3 = \wh{U}_3 (p-\hbar), \;\; p \wh{U}_4 = \wh{U}_4
(p+\hbar) ,
\eeq
\beq{2.13b}
[\wh{U}_1, \wh{U}_2 ] \, = \, [\wh{U}_1, \wh{U}_3 ] \, = \,
[\wh{U}_2, \wh{U}_4 ] \, = \, [\wh{U}_3, \wh{U}_4 ] = 0 ,
\eeq
\beq{2.13c}
[ \wh{U}_1, \wh{U}_4 ] = {\rm Det}U_0 ,\;\;\;\;\; [\wh{U}_3, \wh{U}_2 ] = -
{\rm Det}U_0 ,
\eeq
where ${\rm Det} U_0$ stands for a limit version of (\ref{2.12}):
\beq{2.13d}
 {\rm Det} U_0 = ( \wh{U}_1 \wh{U}_4 - \wh{U}_2 \wh{U}_3 ) \frac{\hbar}{p}
= ( \wh{U}_4 \wh{U}_1 - \wh{U}_3 \wh{U}_2 ) \frac{\hbar}{p} .
\eeq

\begin{prop}\label{P5}
A possible solution for (\ref{2.13a})-(\ref{2.13d}) is 
\beq{2.14}
\wh{U}_1 = \partial_1, \;\; \wh{U}_2 = z_2, \;\; \wh{U}_3 =
-\partial_2, \;\; \wh{U}_4 = z_1;
\eeq
\beq{2.15}
p = \hbar ( z_1 \partial_1 + z_2 \partial_2 +1),
\eeq
where we denote
$\partial_i \equiv \frac{\partial}{\partial z_i}$.
\end{prop}

\rem The representation given by (\ref{2.14})-(\ref{2.15}) is not 
unique. In particular, the rescaling 
$\wh{U}_i\rar c_i \wh{U}_i$ (where $c_i$ are numerical 
constants such that $c_1 c_4 = c_2 c_3$) is allowable.

The Proposition \ref{P5}  together with the connection formula (\ref{2.7}) 
allows us to write out the explicit form of the matrix $U_0$
\beq{2.16}
U_0 = \left(  \ar{cc} \partial_1 & z_2 \\ -\partial_2 & z_1 \er \right)\,
\sqrt{\frac{\hbar}{p}} .
\eeq
Note that this matrix is ``unimodular'', i.e., 
${\rm Det} U_0 = ( \partial_1 z_1 + z_2 \partial_2 )\frac{\hbar}{p} = 1$.

To describe the obtained representation of the algebra ${\cal U}_0$ completely
one has to define a space where operators (\ref{2.14})-(\ref{2.16}) act. It 
is natural to think that this space is $D(z_1,z_2)$ -- a space of 
holomorphic functions of two complex variables.

Let us recall that $D(z_1,z_2)$ is a space spanned on the vectors
\beq{3.2}
|j,m\rangle =\frac{z_1^{j+m}z_2^{j-m}}{\sqrt{(j+m)!(j-m)!}},\;\;\;\;
j=\ar{lllll}0,&\!\!\frac{1}{2},&\!\!1,&\!\!\frac{3}{2},&\!\! ...\er
\;\;\;\;\;m=-j,..,j
\eeq
and equipped with the scalar product
\beq{3.3}
\langle f,g\rangle =\frac{1}{(2\pi i)^2}\int
\overline{f(z_1, z_2)} g(z_1,z_2)
e^{-z_1\bar z_1-z_2\bar z_2}dz_1 d\bar z_1 dz_2 d\bar z_2.
\eeq
The system (\ref{3.2}) is orthonormal with respect to the scalar product
(\ref{3.3}), that is
$ <j,m|j^{\prime},m^{\prime}>=\delta_{j j^{\prime}} \delta_{m m^{\prime}}$.
For the given scalar product a rule of conjugation of operators looks as
follows
\beq{3.5}
(z_i)^{*}=\partial _i,\;\;\;(\partial _i)^{*}=z_i.
\eeq
The question concerning unitarity of the matrix $U_0$ is discussed
in Appendix B.

\subsection*{\red \ns\bf B. Connection with $\bf T^*{\cal B}$. }
\red The generators of $sl(2)$ can be realized on $D(z_1,z_2)$ as
differential operators:
\beq{3.1}
l_{+}=z_1\partial _2,\;\;l_{-}=z_2\partial
_1,\;\;l_3=\ar{c} \frac 12 \er \! (z_1\partial _1-z_2\partial _2).
\eeq
Using these expressions we can compare the representation of the
algebra ${\cal L}$ (or, more precisely, its limit version ${\cal L}_0$)
given by Theorem 1 with the representation given by Proposition 
\ref{P4}.

Indeed, in the limit $\ga\rar 0$ the initial formula (\ref{1.7})
acquires form
\beq{3.05}
L=I+\ga L_0 +O(\ga^2) , \;\;\; L_0 = U_0\left(\ar{cc} p&\\ &-p\er\right)
U_0^{-1} .
\eeq

Substituting here the explicit expressions (\ref{2.15})-(\ref{2.16})
for $p$, $U_0$ and using the representation
(\ref{3.1}) for generators of $sl(2)$, one derives the following
limit form of the $L$-operator:
\beq{3.15}
L_0=\hbar \ \left(\ar{cc} 2 + z_1\,\partial_1 - z_2\,\partial_2 &
 2 z_2\,\partial_1 \\  2 z_1\,\partial_2 & 2 - z_1\,\partial_1 + z_2\,
\partial_2 \er \right) =
2\hbar \ \left(\ar{cc} 1+l_3 & l_{-} \\ l_{+} & 1-l_3 \er \right) .
\eeq
Notice that (\ref{3.15}) exactly coincides with
(\ref{L5}) taken in the limit $\ga \rar 0$. This explains why we had to
introduce the factor $q^2$ in (\ref{L5}).

The next observation concerning the limit of $L$-operator reads as
follows.
\begin{prop}\label{P6}
The matrix $L_0$ in the representation (\ref{3.15}) admits
the decomposition
\beq{3.35a}
   L_0 = A_0\,B_0\,A_0^{-1},
\eeq
where
\beq{3.35b}
  A_0 = \left(\ar{cc} z_1^{-1/2} & -z_1^{-1/2} z_2 \\
    0 & z_1^{1/2} \er\right), \;\;\;
 B_0= \hbar\left(\ar{cc} p/\hbar+1/2 & 0 \\ 2\,\partial_2 &
    -(p/\hbar-1/2)\er\right)
\eeq
and $p$ is defined as in (\ref{2.15}).
\end{prop}
This statement can be verified directly.

Let us comment on the meaning of this proposition. First, note that $A_0$
is a realization of a group-like element of the Borel subgroup of $SL(2)$.
Moreover, this explicit form of $A_0$ is straightly connected with the
construction of the model space $\cal M$ developed by Gelfand et al.\REF{1}
Indeed, the space $D(z_1,z_2)$ being a realization of the model
space for $SL(2)$ (compare (\ref{0.0}) and (\ref{3.2})) is spanned on
monomials with arguments which are combinations of the entries of $A_0$.
On the other hand, $B_0$ is of opposite (with respect to $A_0$)
triangularity and its entries are operators acting on a given realization
of the model space. Therefore, $B_0$ can be regarded as an element of the
space dual to the corresponding Borel subalgebra.

Thus, $A_0$ and $B_0$ are coordinates in the base and in a fiber of the 
cotangent bundle $T^*\B$.
At this stage the appearance of $T^*\B$ "inside" the algebra $\cal L$ looks
somewhat mysterious, but we shall clarify it later.

\subsection*{\red \ns\bf C. Clebsch-Gordan coefficients. }
\red  Let us consider an action of the generators of the algebra
${\cal U}_0$ defined in (\ref{2.15})-(\ref{2.16}) on the space $D(z_1,z_2)$
(which is a realization of the model space).
The action of these operators on the basic vectors (\ref{3.2}) is given by
\beq{3.8}
p \; |j,m\rangle =(2j+1)\hbar \; |j,m\rangle ,
\eeq
\beq{3.9} \!\!\!\! \ar{ll}
U_1|j,m\rangle =\left(\frac{j+m}{2j+1}\right)^{1/2}|j-\frac 12 ,m-\frac 12
\rangle , &
U_2|j,m\rangle =\left(\frac{j-m+1}{2j+1}\right)^{1/2}|j+\frac 12 ,m-\frac 12
\rangle ,\\     \\
U_3|j,m\rangle =-\left(\frac{j-m}{2j+1}\right)^{1/2}|j-\frac 12 ,m+\frac 12
\rangle , &
U_4|j,m\rangle =\left(\frac{j+m+1}{2j+1}\right)^{1/2}|j+\frac 12 ,m+\frac 12
\rangle .  \er
\eeq

%%%%%%%%%

Formula (\ref{3.8}) allows us to identify the operator $p$ as
$p=2\wh{j}+1$, where $\wh{j}$ is the operator of spin. Hence,
invariant subspaces of $p$ on the model space are those with
fixed value of spin $j$.

Formulae (\ref{3.9}) show that $U_i$ are generators of the basic shifts
on the model space (as illustrated on Fig.1). This observation is very
important. As we shall see later, the same picture holds for $q\neq 1$.

%%%%%%%%%%%%%%%%   Fig.1.  %%%%%%%%%%%%%%%%%%%%%%%
\vspace*{1mm}
\hbox{
\begin{picture}(200,150)(-100,0)
\put(20,20){\vector(1,0){150}} \put(70,20){\vector(0,1){131}}

\multiput(90,20)(0,10){12}{\line(0,1){5}}
\multiput(110,20)(0,10){12}{\line(0,1){5}}
\multiput(130,20)(0,10){12}{\line(0,1){5}}

\multiput(70,80)(10,0){8}{\line(1,0){5}}
\multiput(70,100)(10,0){8}{\line(1,0){5}}
\multiput(70,120)(10,0){8}{\line(1,0){5}}

\put(73,147){$j$} \put(173,23){$m$}

\put(72,5){$m\!-\!\frac{1}{2}$} \put(106,5){$m$}
\put(122,5){$m\!+\!\frac{1}{2}$}

\put(42,78){$j\!-\!\frac{1}{2}$} \put(42,98){$j$}
\put(42,118){$j\!+\!\frac{1}{2}$}

\put(75,65){$U_1$} \put(75,125){$U_2$}
\put(135,65){$U_3$} \put(135,125){$U_4$}

\thicklines
\put(110,100){\vector(1,1){20}} \put(110,100){\vector(-1,1){20}}
\put(110,100){\vector(1,-1){20}} \put(110,100){\vector(-1,-1){20}}

\end{picture} } \vspace*{0.5mm}
\begin{center}
\parbox{9cm}{\small {\bf Fig.1:} \
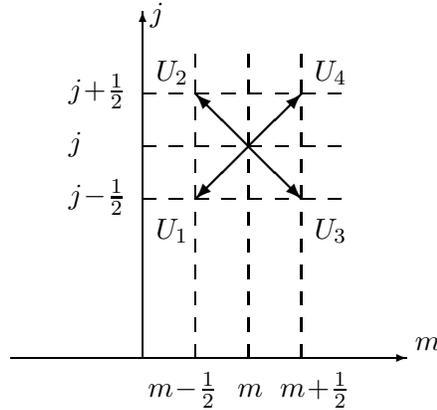
 Action of the operators $U_i$ on the model space.}
\end{center}
%%%%%%%%%%%%%%%%% End of Fig.1. %%%%%%%%%%%%%%%%%%%

Now comparing the matrix elements $\langle\,j'',m''|U_i|j,m\,\rangle$
following from (\ref{3.9})
with values of the Clebsch-Gordan coefficients (CGC) for decomposition of
the tensor product of irreducible representations $V_j\otimes V_{\frac 12}$
for $sl(2)$ which are given by the Van-der-Waerden formula
\beq{3.10} 
\left\{\ar{ccc} j & \frac 12 & j^{\prime \prime } \\ m &
m^{\prime } & m^{\prime \prime } \er \right\} = \delta_{m^{\prime \prime },
m+m^{\prime }} \sqrt{ \frac{(j+\frac 12-j^{\prime \prime })!
(j+j^{\prime \prime }-\frac 12 )!(j^{\prime \prime }+\frac 12-j)! }
{(j+j^{\prime \prime }+\frac 32)!} } \; \times
\eeq
\[ \!\! \times \sum\limits_{r\geq 0}
\frac{(-1)^r \sqrt{ (j+m)!(j-m)!(j^{\prime \prime}+m^{\prime \prime})!
(j^{\prime\prime }-m^{\prime \prime })!(2j^{\prime\prime }+1) } }
{r!(j+\frac 12-j^{\prime \prime}-r)!(j-m-r)!(\frac 12+m^{\prime }-r)!
(j^{\prime \prime }-\frac 12+m+r)!(j^{\prime \prime }-j-m^{\prime }+r)!},  \]
we establish the following correspondence
\beq{3.14} \! \ar{ll}
\langle j^{\prime \prime} , m^{\prime \prime} \mid U_1\mid j, m\rangle=
\delta_{j^{\prime \prime },j- \frac 12 } \left\{\ar{ccc} j & \frac 12 &
j^{\prime \prime} \\ m & -\frac 12 & m^{\prime \prime} \er \right\} , \\ [1mm]
\langle j^{\prime \prime} , m^{\prime \prime} \mid U_2 \mid j, m\rangle=
\delta_{j^{\prime \prime },j+ \frac 12 }\left\{\ar{ccc} j & \frac 12 &
j^{\prime \prime} \\ m & -\frac 12 & m^{\prime \prime} \er \right\} , \\ [1mm]
\langle j^{\prime \prime} , m^{\prime \prime} \mid U_3 \mid j, m\rangle=
\delta_{j^{\prime \prime },j- \frac 12 }\left\{\ar{ccc} j & \frac 12 &
j^{\prime \prime} \\ m & \frac 12 & m^{\prime \prime} \er \right\}, \\ [1mm]
\langle j^{\prime \prime} , m^{\prime \prime} \mid U_4 \mid j, m\rangle=
\delta_{j^{\prime \prime },j+ \frac 12 }\left\{\ar{ccc} j & \frac 12 &
j^{\prime \prime} \\ m & \frac 12 & m^{\prime \prime}\er \right\}. \er
\eeq

Thus, we proved the following statement:
\begin{prop}\label{P7}
The generators $U_i$ of the algebra ${\cal U}_0$ are operators of the 
basic shifts on the model space for $sl(2)$ and they generate the 
Clebsch-Gordan coefficients corresponding to decomposition of the product 
$V_j\otimes V_{\frac 12}$ of the irreps of sl(2).
\end{prop}

This statement allows to call the matrix $U_0$ a ''generating matrix'' 
(by analogy with the notion of a generating function) for CGC. 

\rem Usually, introducing a generating object (well-known examples are
the generating functions for different sets of polynomials, e.g., for the 
Legendre polynomials), one makes properties of the objects under 
consideration more evident. We think that the notion of generating 
matrix will be useful for calculations involving CGC of classical and 
quantum algebras.
 
\subsection*{\red \ns\bf D. Wigner-Eckart theorem.}

\red One should underline a connection of the results obtained above
(Proposition \ref{P7}) and the well-known mathematical construction --
Wigner-Eckart theorem,\REF{19} which has important applications 
in quantum mechanics.

Let us remind that the Wigner-Eckart theorem gives CGC for classical
Lie algebra ${\cal J}$ as matrix elements of some set of operators. These 
operators are called {\it tensor operators}. They map the corresponding 
model space $\cal M$ onto itself and have special transformation properties
under adjoint action of the algebra. In the case of ${\cal J}=sl(2)$ the 
Wigner-Eckart theorem reads as follows.
\begin{theorem} \label{Twe}
Let $l_{+}$, $l_{-}$ and $l_{3}$ be the generators of $sl(2)$ and
let $T^{j\;}_{\;m}$, $m=-j,..,j$ be a set of operators acting on $\cal M$
and obeying the commutation relations
\beq{we1}
 [l_3,T^{j\;}_{\;m}]=m\,T^{j\;}_{\;m},\;\;\;\;
 [l_{\pm},T^{j\;}_{\;m}]=\sqrt{(j\mp m)(j\pm m+1)}\; T^{j\;}_{\;m\pm 1},
\eeq
where $j(j+1)$ is an eigenvalue of the Casimir operator for $sl(2)$.
Then the matrix elements of $T^{j\;}_{\;m}$ on $\cal M$ are proportional to
Clebsch-Gordan coefficients:
\[ <j^{\prime\prime} m^{\prime\prime} |T^{j\;}_{\;m}|j^{\prime} m^{\prime} >
= C_{j^{\prime\prime}j^{\prime}}^{j}\;\left\{\ar{ccc}j^{\prime} & j &
j^{\prime\prime}\\ m^{\prime} & m & m^{\prime\prime} \er \right\},    \]
where the coefficients $C_{j^{\prime\prime}j^{\prime}}^j$ do not depend on
$m$, $m^{\prime}$, $m^{\prime\prime}$.
\end{theorem}

Proposition \ref{P7} says that any tensor operators of spin
$j=1/2$ (that is $\{T^{1/2}_{1/2},\, T^{1/2}_{-1/2}\}$,\ 
$T^{1/2}_m\, :\, V_j \mapsto V_j\otimes V_{1/2}=V_{j+1/2}\oplus V_{j-1/2}$)
may be constructed via the operators $U_i$ (in fact, it is evident from 
Fig.1). Indeed, comparing the commutation relations obtained directly from
(\ref{2.16}) and (\ref{3.1})
\beq{we20} \ar{llll}
{}[l_{+},U_1]=U_3, & [l_{+},U_2]=U_4, & [l_{+},U_3]=0, & [l_{+},U_4]=0,\\
& & & \\ {} [l_{-},U_1]=0, & [l_{-},U_2]=0, & [l_{-},U_3]=U_1, & [l_{-},U_4]
= U_2,\\ & & & \\ {}[l_{3},U_1]=-\frac{1}{2} U_1, &[l_{3},U_2] =-\frac{1}{2}
U_2, & [l_{3},U_3] = \frac{1}{2} U_3, & [l_{3},U_4]= \frac{1}{2} U_4  \er
\eeq
with Theorem \ref{Twe}, we get the following.
\begin{prop}\label{P8}
The generators $U_i$ of the algebra ${\cal U}_0$ form a basis for tensor 
operators of spin $1/2$, that is components $T^{1/2}_{1/2}$ and 
$T^{1/2}_{-1/2}$ of any tensor operator of spin 1/2 can be realized as
linear combinations of $U_i$:
\beq{3.55}
T^{1/2}_{-1/2}= \mu(p) \;U_1+\nu(p) \;U_2,\;\;\; 
T^{1/2}_{1/2}=\mu(p) \;U_3+\nu(p) \;U_4,
\eeq
where $\mu(p)$ and $\nu(p)$ are functions only of $p=2\wh{j}+1$.
\end{prop}

\section{ \ns\bf DEFORMED CASE. }
\setcounter{equation}{0}

\red Now we want to extend the results
obtained in the previous section to the case of $q\neq 1$. In
particular, we are going to examine the representations of the algebra
${\cal U}$ (see Definition \ref{DU} above) and to show that the
corresponding matrix $U$ generates Clebsch-Gordan coefficients for the
deformed Lie algebra. For these purposes we shall exploit a natural
connection of ${\cal U}$ with $(T^*\B)_q$.

\subsection*{\red \ns\bf A. The $\bf q$-oscillators approach. }
\red There exist different ways to obtain desirable representations of
the algebra ${\cal U}$. First we describe a more direct but less
instructive method, which is similar to that used in the
non-deformed case.

By analogy with the non-deformed case studied above, one can assume that
the entries of the matrix $U$ might be realized as operators
(deformations of those obtained in Proposition \ref{P5}) acting on
the space of two complex variables. Indeed, using the definition
(\ref{2.6}) of the central element of $\cal U$ and taking into account
the identity for $q$-numbers
\beq{q1} [a]\,q^b+[b]\,q^{-a}=[a+b] ,
\eeq we can rewrite (\ref{2.8})-(\ref{2.11}) in the following way:
\beq{q2}
\wh{U}_1 \wh{U}_3 = q^{-1} \wh{U}_3 \wh{U}_1, \;\;\, \wh{U}_2 \wh{U}_4
= q^{-1} \wh{U}_4 \wh{U}_2, \;\;\, \wh{U}_1 \wh{U}_2 = \wh{U}_2
\wh{U}_1, \;\;\, \wh{U}_3 \wh{U}_4 = \wh{U}_4 \wh{U}_3 , \eeq \beq{q3}
\wh{U}_1 \wh{U}_4 - q^{-1} \wh{U}_4 \wh{U}_1 = q^{-1} {\rm Det}U \,
q^{p/\hbar},\;\;\;\;\; \wh{U}_3 \wh{U}_2 - q\wh{U}_2 \wh{U}_3 =
 - {\rm Det}U \,q^{-p/\hbar}.  
\eeq

The relations (\ref{q3}) are well known in the theory of
$q$-oscillators ($q$-bosons).\REF{20} Recall that $q$-analogues of
creation, annihilation, and number operators form the deformed
Heisenberg algebra defined by the commutation relations
\beq{q3.1}
a\, a^+ - q \, a^+ a = N^{-1},\;\; N\, a = q^{-1}\, a\, N,\;\;
N\, a^+ = q\, a^+ N \, ,
\eeq
and they can be realized in terms of multiplication and
difference operators:
\beq{q4}
a^+=z,\;\;a=z^{-1} [z\,\partial_z],\;\; N = q^{z\,\partial_z} .  
\eeq

Using two pairs of generators of the deformed Heisenberg algebra, one
can construct the generators of $U_q(sl(2))$: 
$l_+ = a_1^+\,a_2$, $l_- = a_2^+\,a_1$, $q^{l_3}=N_1^{1/2}N_2^{-1/2}$. 
Applying here the representation (\ref{q4}) one gets 
\beq{q33} 
l_+ = z_1 z_2^{-1} [z_2\partial_2 ],\;\;\;\;\; 
l_- = z_2 z_1^{-1} [z_1\partial_1 ],\;\;\;\;\;
q^{l_3} = q^{\frac{1}{2} (z_1 \partial_1 - z_2 \partial_2)}.  
\eeq 
The Casimir operator (\ref{L8}) of $U_q(sl(2))$ in this realization is 
given by 
\beq{q4.5}
C= q\,N_1 N_2 + q^{-1} N_1^{-1} N_2^{-1}.  
\eeq

Now, comparing, (\ref{q2})-(\ref{q3}) with (\ref{q3.1}), it is easy to
conclude that the pairs $(\wh{U}_1,\wh{U}_4)$ and
$(\wh{U}_2,\wh{U}_3)$ are similar to two pairs of $q$-boson operators.

Taking into account the Weyl-like form of relations (\ref{q2}) and
having already found explicit expressions (\ref{2.14})-(\ref{2.15})
for the generators of algebra ${\cal U}_0$, we get an answer for $D$
and $\wh{U}$ in terms of $q$-oscillators. More precisely, a
straightforward calculation allows to verify the following statement:
\begin{prop}\label{P9}
Equations (\ref{q2})-(\ref{q3}) have the family of solutions: 
\beq{q6}
q^{p/\hbar} = q\,N_1\,N_2,\;\;\;\; \wh{U} = \left( \ar{cc} \al_0 \,a_1
\, N_1^{\al}\,N_2^{-\be} &\be_0\, a_2^+\,N_1^{\be}\,N_2^{-\al} \\ [1mm]
-\ga_0\, a_2 \,N_1^{-(1+\be)}\, N_2^{\al} &
\dl_0\,a_1^+\,N_1^{-\al}\,N_2^{1+\be} \er\right), 
\eeq 
where $\;\al_0\,\dl_0=q\,\be_0\,\ga_0$.
\end{prop}
Let us note that this form of $\wh{U}$ is consistent with the
condition (\ref{2.5}).

Taking into account the connection formula (\ref{2.7}) and applying to
the generators $a_i$, $a_i^+$, $N_i$ the representation (\ref{q4}),
one obtains from (\ref{q6}) a family of representations of the algebra
${\cal U}$. To select some of them, we have to impose an additional
condition.

 As mentioned above (see (\ref{3.05})-(\ref{3.15})), in the non-deformed
case substitution of the generating matrix $U_0$ in the formula
(\ref{1.7}) gives the matrix $L_0$ which exactly coincides with the
limit version of the matrix (\ref{L5}). It is natural to suppose that
the generating matrix $U$ corresponding to deformed algebra produces
in the same way the matrix (\ref{L5}) itself. Bearing in mind the
property (\ref{2.130}), we obtain the following
\begin{prop}\label{P10}
The condition $\wh{U}D\wh{U}^{-1}=L$, where $L$ is the matrix
(\ref{L5}), $D$ is given by 
\beq{q5.5} 
 D = \left(\ar{cc} q^{p/\hbar} & \\ & q^{-p/\hbar} \er\right) =
 \left(\ar{cc} q\,N_1\,N_2 & \\ &
 q^{-1}\,N_1^{-1}\,N_2^{-1} \er\right) \, ,
\eeq and $\wh{U}$ is given by (\ref{q6}), imposes the following
restrictions:
\beq{q8}
\al+\be+\frac{1}{2}=0,\;\;\; \al_0=q\,\ga_0 ,\;\;\; \be_0= \dl_0.
\eeq
\end{prop}
Substitution of (\ref{q8}) into (\ref{q6}) completes a description of
$\wh{U}$ in terms of $q$-oscillators.

\subsection*{\red \ns\bf B. Connection with $\bf (T^*{\cal B})_q$. }
\red Now we are going to develop another approach to constructing
representations of $\cal U$. It is more universal since it is based on
the connection (which takes place for arbitrary quantum Lie algebra)
of the algebra ${\cal L}$ (see Definition \ref{DL}) with $(T^*\B)_q$
and on the interpretation of the deformed Borel subgroup $\B_q$ as a
quantum model space.

To clarify the announced connection we start with the following
theorem (this is a version of the theorem given in Ref.10 for
$L$-operators with nonultralocal relations)
\begin{theorem} \label{ABL}
 Let the matrices $A$ and $B$ obey the relations of type (\ref{1.1}):
\beq{4.1}
R_{\pm }\up{1}{A} \, \up{2}{A} \,=\, \up{2}{A} \, \up{1}{A}
R_{\pm} \ , \;\;\;\;
R_{\pm } \up{1}{B} \, \up{2}{B} \,=\, \up{2}{B} \, \up{1}{B} R_{\pm}
\eeq
and the additional exchange relation
\beq{4.2}
\up{1}{A} \, \up{2}{B} \,=\, \up{2}{B} \, \up{1}{A} R_{+} \ , \;\;\;\;
\up{2}{A} \, \up{1}{B}R_{-} \,=\, \up{1}{B} \, \up{2}{A} \  .
\eeq

Then the $L$-operator constructed by means of similarity transformation
\beq{4.3} 
L=ABA^{-1} 
\eeq 
satisfies the relation (\ref{1.3}).
\end{theorem}
\rem Since (\ref{4.1}) defines a quantum group structure, $A^{-1}$ in
(\ref{4.3}) should be understood as an antipode of $A$.

\proof of Theorem \ref{ABL} is straightforward :
\[ \ar{c}
\up{1}{L} R^{-1}_{-} \up{2}{L}R_{-}= \up{1}{A}
\up{1}{B}(\up{1}{A})^{-1} R_{-}^{-1} \up{2}{A} \up{2}{B}
(\up{2}{A})^{-1} R_{-}= \up{1}{A} \up{1}{B} \up{2}{A} R_{-}^{-1}
(\up{1}{A})^{-1} \up{2}{B} (\up{2}{A})^{-1}R_{-}=\\ \\
=\up{1}{A}\,\up{2}{A}\,\up{1}{B}\,\up{2}{B}R_{+}^{-1}(\up{1}{A})^{-1}
(\up{2}{A})^{-1} R_{-}= R_{+}^{-1}
\up{2}{A}\,\up{1}{A}\,\up{2}{B}\,\up{1}{B} R_{-} (\up{2}{A})^{-1}
(\up{1}{A})^{-1}= \\ \\ = R_{+}^{-1} \up{2}{A} \up{2}{B}\up{1}{A}R_{+}
(\up{2}{A})^{-1} \up{1}{B} (\up{1}{A})^{-1}= R_{+}^{-1}\up{2}{A}
\up{2}{B} (\up{2}{A})^{-1} R_{+} \up{1}{A} \up{1}{B}(\up{1}{A})^{-1}=
R_{+}^{-1} \up{2}{L} R_{+}\up{1}{L}.\er \]

Thus, for a given quantum group $G_q$, the algebra $\cal L$ is
embedded into the algebra generated by entries of $A$ and $B$ obeying
(\ref{4.1})-(\ref{4.2}).  To argue that (\ref{4.1})-(\ref{4.2})
describe a $q$-analogue of $T^*\B$, let us notice that the
non-symmetric (with respect to $R$-matrices) form of the relations
(\ref{4.2}) imposes some restriction on the structure of the matrices
$A$ and $B$. Say, if $R_+$ is an upper triangular matrix, then $A$ and
$B$ must be upper and lower triangular, respectively. Therefore, one
may think of $A$ and $B$ as coordinates in the deformed Borel
subgroup ${\cal B}_q$ and in the dual quantum space, respectively. In
other words, the matrices $A$ and $B$ are coordinate and momentum on
the deformed phase space $(T^*\B)_q$ respectively. Thus
(\ref{4.1})-(\ref{4.2}) may be regarded as a definition of
$(T^*\B)_q$ (for additional comments see Ref.10).

We should underline here that, although the matrices $A$ and $B$ look
similarly on quantum level, they
transform into different objects when $q\rar 1$. Indeed, in the limit
$q\rar 1$ one has $L \rar I + \ga\hbar L_0$ and the corresponding
limit forms of $A$ and $B$ are
\beq{4.45}
 A \rar A_0, \;\;\; B \rar I+\ga\hbar B_0,
\eeq
where $A_0$ is a group-like element, whereas
$B_0$ is rather an element of algebra (see (\ref{3.35b}) as an example
of $A_0$, $B_0$ for $sl(2)$).

Comparing the statements of Theorems 1 and \ref{ABL} and
taking into account the equality (\ref{2.130}), we get the formula
\beq{4.115}
 L = A\, B \, A^{-1} = \wh{U} \, D \, \wh{U}^{-1},
\eeq
which points out a possibility to construct the matrix $\wh{U}$
obeying (\ref{2.8})-(\ref{2.9}) via the generators of $(T^*\B)_q$. 
This connection is very important; below we consider it for 
$SL_q(2)$ in all details.

Now let us turn to the example of $SL_q(2)$. For $R_{\pm}$ defined as
in (\ref{1.4}) one can choose \beq{4.4} A=\left(\ar{cc} a & c \\ 0 &
a^{-1} \er \right) ,\ \ \ \ B=\left(\ar{cc} b & 0 \\ d & b^{-1}
\er\right) .  \eeq

Explicit relations for the generators of $(T^*\B)_q$ following from
(\ref{4.1})-(\ref{4.2}) are
\beq{4.5} a\,c=q^{-1}c\,a,\ \ \ b\,c=q^{1/2}c\,b,\ \ \ a\,b=q^{1/2}b\,a;
\eeq
\beq{4.6} b\,d=q^{-1}d\,b,\ \ \ a\,d=q^{1/2}d\,a,\ \ \
c\,d=q^{-1/2}d\,c+q^{-1/2}\omega\, b^{-1}a.
\eeq

Performing the following decomposition
\beq{4.7}
d=d_0+d_1=d_0+q^{1/2}c^{-1}b^{-1}a \, ,
\eeq
we transform (\ref{4.6}) to homogeneous form:
\beq{4.8} b\,d_0=q^{-1}d_0\, b,\ \ \ \ a\, d_0=q^{1/2}d_0 \, a,
 \ \ \ \ c\, d_0=q^{-1/2}d_0 \, c .
\eeq

Thus, (\ref{4.5}) and (\ref{4.8}) describe four variables obeying
Weyl-like commutation relations. Using the jargon of conformal field
theory, we shall call these formulae "free field representation" and
the generators $a$, $b$, $c$, $d_0$ "free field" variables.

\rem The last of equations (\ref{4.6}) is nothing but a commutation
relation entering the definition of deformed Heisenberg algebra.
Indeed, comparing (\ref{4.5})-(\ref{4.6}) with (\ref{q3.1}), one can
establish the following correspondence ($\rho$ stands for arbitrary
numerical constant):
$$
 c \sim N^{\rho}\,a^{+},\;\;\; d \sim -\omega \,
 N^{-1/2-\rho}\,a,\;\;\; b^{-1}\,a \sim q^{\rho}\,N^{-3/2}.
$$
Thus, the transformation (\ref{4.7}) can be interpreted as
"bosonization" of $q$-oscillators.

Now, substituting (\ref{4.4}) in (\ref{4.3}), we get 
\beq{4.9} \ar{c}
L = q^{1/2} \left( \ar{cc} a & c \\ 0 & a^{-1} \er \right) \left(
\ar{cc} b & 0 \\ d & b^{-1} \er \right) \left( \ar{cc} a^{-1} & -q\,c
\\ 0 & a \er \right) = \\ \\ = \left( \ar{cc} q\, (b+b^{-1}) +
a^{-1}c\, d_0 & - q^2 a\, c\, (b+q \, a^{-1} c\, d_0) \\ (a\,c)^{-1}
(b^{-1} + q^{-1} a^{-1}c\, d_0) & -q^2 a^{-1} c\, d_0 \er \right). \er
\eeq 
This matrix provides a "free field" realization of the algebra
$\cal L$ for $U_q(sl(2))$. Note that the additional scaling factor
$q^{1/2}$ was introduced in (\ref{4.9}) to ensure a coincidence of the
Casimir operators calculated by formulae (\ref{1.6}) for the matrix
(\ref{4.9}): 
\beq{4.10} 
K_1 = q^2 (b+b^{-1}), \;\;\; K_2 = q^3  
\eeq
with those for the matrix (\ref{L5}). In fact, we redefined the matrix
$B$ in (\ref{4.4}) as \beq{4.95} \wt{B} = q^{1/2} B .  \eeq Comparing
the Casimir operator $K_1$ given by (\ref{4.10}) with one given by
(\ref{L6}), we identify the operator $b$ with the power of the operator
of spin $\wh{j}$ :
\beq{4.105} 
b = q^{2 \wh{j} +1}.  
\eeq

It follows from (\ref{4.10}) that matrix $L$ contains only three
independent variables (it is easy to see from the explicit form
(\ref{4.9}) that these are $b$, $ac$ and $a^{-1}cd_0$).  Moreover,
direct calculation using (\ref{4.5}),(\ref{4.8}) shows that all
elements of the matrix $L$ commute with operator $b$. That agrees with
the property (\ref{1.12}).

Now exploiting the connection described by formula (\ref{4.115}), one
can obtain  an exact expression for $\wh{U}$.
\begin{theorem} \label{main}
The algebra ${\cal U}\equiv \{\wh{U},p\}$ with defining relations
(\ref{2.8})-(\ref{2.9}) has the following realization in terms of
generators $a$, $b$, $c$, $d_0$ 
\beq{5.00} 
b = q^{p/\hbar},\;\;\;\;
\wh{U} = \left( \ar{cc} \frac{1}{\omega}\, a\,(b+a^{-1}\,c\,d_0)\,
e^{-\frac{i \xi}{2}} & c\,e^{\frac{i \xi}{2}} \\ [1mm]
\frac{1}{\omega}\, c^{-1}(b^{-1} + q^{-1} a^{-1}\,c\,d_0)\,
e^{-\frac{i\xi}{2}} & a^{-1}\, e^{\frac{i \xi}{2}} \er \right), 
\eeq 
where $\omega\equiv q-q^{-1}$,\ $d_0$ is defined in (\ref{4.7}), and
\beq{5.15}
e^{i \xi}= a^{-1}\;b^{\ga}\;c^{-1}\;d_0^{-1}
\eeq 
with $\ga$ being an arbitrary constant.
\end{theorem}

This theorem gives a "free field" representation of the algebra $\cal
U$.  Let us remark that the remaining freedom in (\ref{5.15})
corresponds only to canonical transformations (since $\xi$ and $p$ are
conjugate variables).

The formulated theorem will be proved in several steps. First, we
introduce a lower-triangular matrix which diagonalizes the matrix
$\wt{B}$: 
\beq{5.1} 
V = \left( \ar{cc} v_1 & 0 \\ v_3 & v_2 \er \right), \;\;\; 
\wt{B} = V \wt{B}_0 V^{-1}, \;\;\; \wt{B}_0 = \left(\ar{cc} q^{1/2} b 
& 0 \\ 0 & q^{1/2} b^{-1} \er \right) \equiv q^{1/2} B_0 .  
\eeq
\begin{prop}
A possible solution for the matrix $V$ is \beq{5.3} v_1 =v_1(b),\;\;\;
v_2 =v_2(b),\;\;\; v_3 = d \,v_1(b)\,f(b), \eeq where $v_1(b)$,
$v_2(b)$ are arbitrary functions of $b$ and 
$f(b) = (b-qb^{-1})^{-1}$.
\end{prop}

Thus, matrix $L$ given by (\ref{4.9}) admits a decomposition of the form:
\beq{5.5} L=\wh{U}_0\, \wt{B}_0\, \wh{U}_0^{-1}, \;\;\;\; \wh{U}_0 =
A\,V.  \eeq

However, this diagonalization is not unique. Using an arbitrary power
of the diagonal matrix $Q$, which depends on the variable conjugate to
$b$, 
\beq{5.6} 
Q = \left( \ar{cc} e^{i\xi} & \\ & e^{-i\xi} \er \right), 
\;\;\;\; b\,e^{i\xi} = q\,e^{i\xi}\,b, 
\eeq 
we obtain a
family of diagonalizing matrices: 
\beq{5.7}
L= \wh{U}_{\dl}\,\wt{B}_{\dl}\,\wh{U}_{\dl}^{-1}, \;\;\; \wh{U}_{\dl}=
A\,V\,Q^{\dl}, \;\;\; \wt{B}_{\dl} = Q^{-\dl} \wt{B}_0 \, Q^{\dl} =
q^{\dl}\,\wt{B}_0 = q^{\dl+1/2}\, B_0.  
\eeq 
An explicit form of the diagonalizing matrix is 
\beq{5.8} 
\wh{U}_{\dl} = A\,V\,Q^{\dl} = \left( \ar{cc} (a\,v_1 +c\,d\,v_1 f)\, 
e^{i\dl \xi} & c\,v_2 \,e^{-i\dl \xi} \\ a^{-1} d\,v_1 f\, e^{i\dl \xi} & 
a^{-1}v_2 \, e^{-i\dl \xi} \er \right).  
\eeq

Here we should describe a new object $e^{i \xi}$ which appeared in the
matrix $\wh{U}$. We assume that the following Weyl-like relations
hold: 
\beq{5.10} 
a\, e^{i \xi}= q^{\al} e^{i \xi}\, a, \;\;\; b\, e^{i
\xi}= q\,e^{i \xi}\, b, \;\;\; c\, e^{i \xi}= q^{\be} e^{i \xi}\,
c,\;\;\; d_0\, e^{i \xi}= q^{\ga} e^{i \xi}\,d_0.  
\eeq
\begin{prop}
The set of equations (\ref{5.10}) is equivalent to
\beq{5.11}
e^{i\xi}= a^{\be+(\ga-1)/2}\;b^{\ga}\;c^{(\ga-1)/2-\al}\;d_0^{-1}.
\eeq
\end{prop}

Now we have to remind that the matrix $U$ (and $\wh{U}$ as well)
described in the Theorem 1 has to satisfy the relation (\ref{1.8}) or,
equivalently, to the relation
\beq{5.85}
 \up{1}{B_0} \,\up{2}{\wh{U}_{\dl}} =
  \up{2}{\wh{U}_{\dl}} \, \up{1}{B_0} \, \sigma,
\eeq
where $\sigma$ and $B_0$ were introduced in (\ref{1.8})
and (\ref{5.1}), respectively. A straightforward calculation using
(\ref{4.5})-(\ref{4.6}) leads to the following.
\begin{prop}
The matrix $\wh{U}_{\dl}$ given by (\ref{5.8}) satisfies the relation
(\ref{5.85}) only for $ \dl = -1/2$.
\end{prop}
It is worth mentioning that such a choice of $\dl$ exactly compensates
the renormalization of the matrix $B$ in (\ref{4.95}), i.e., \
$\wt{B}_{-1/2} = B_0$.

Bearing in mind the formula (\ref{4.7}), one can
rewrite (\ref{5.8}) for $\dl = \ar{c}-\frac{1}{2} \er$ as follows
\beq{5.9} 
\wh{U} \equiv \wh{U}_{-1/2} = \left( \ar{cc}
a(b+a^{-1}\,c\,d_0)\,w\, e^{-\frac{i \xi}{2}} & c\,v\,e^{\frac{i
\xi}{2}} \\ c^{-1}(b^{-1} + q^{-1} a^{-1}\,c\,d_0)\,w\,e^{-\frac{i
\xi}{2}} & a^{-1}\,v \, e^{\frac{i \xi}{2}} \er \right), 
\eeq 
where $w \equiv f(b) v_1(b)$,\  $v \equiv v_2(b)$.

Finally, a direct check shows (see the Appendix C) that the matrix
(\ref{5.9}) obeys eqs.$\!$ (\ref{2.8})-(\ref{2.9}) if the functions
$w$, $v$ are constant (we chose them as follows: $v(b)=1$,
$w(b)=\frac{1}{\omega}$) and the coefficients in
(\ref{5.10})-(\ref{5.11}) satisfy the conditions $\be \!=\! -\al$, \
$\ga \!=\! \al-\be-1 \!=\! 2\al-1$.  Thus, Theorem \ref{main} is
proven.

Let us end the discussion of relation of $(T^*\B)_q$ to algebras
$\cal L$ and $\cal U$ with one more statement:
\begin{theorem}
The algebra generated by coordinates on $(T^*\B)_q$ is isomorphic to
the algebra generated by entries of the matrix $L$ and $Q$.
\end{theorem}
\proof. Formulae (\ref{4.9}) and (\ref{5.15}) provide explicit
expressions for entries of $L$ and $Q$ via the generators $a$, $b$, $c$,
$d_0$ (up to unessential canonical transformation in (\ref{5.15})).
Conversely, suppose matrix $L$ and the element $e^{i\xi}$ are given.
Then, as it follows from (\ref{4.9}), one can construct from entries of
$L$ the combinations $b$, $ac$ and $a^{-1}cd_0$. Together with (\ref{5.15})
this allows to recover the ``coordinates'' $a$, $b$, $c$, $d_0$.

Although we considered this theorem only for the case of $SL_q(2)$, there
is an evidence that it holds for the generic case. Indeed, in the case
of $G_q=SL_q(N)$ a point on the quantum bundle $(T^*\B)_q$ is parameterized
by $N\times N$ matrices $A$ and $B$. As above, the matrix $L=ABA^{-1}$
satisfies (\ref{1.3}) and therefore its entries generate the corresponding
algebra $\cal L$. However, the dimension of $(T^*\B)_q$ exceeds the
dimension of $\cal L$ : \
${\rm dim}\, (T^*\B)_q - {\rm dim}\, {\cal L} =
%% 2\left( N(N+1)/2 -1 \right) - (N^2-1) =
(N^2+N-2)-(N^2-1)=N-1$.
It us very probable that the remaining $(N-1)$ generators are exactly
those that enter the diagonal unimodular $N\times N$ matrix $Q$.

\subsection*{\red \ns\bf C. Explicit representation.}
\red Now we face the problem of constructing of an explicit
representation for the generators $a$, $b$, $c$, $d_0$. A Weyl-like
form of the commutation relations (\ref{4.5}),(\ref{4.8}) points out
the possibility of getting a realization for these generators in terms
of two pairs of canonical variables.  This also means (due to the
interpretation of (\ref{4.7}) as "bosonization" of $q$-oscillators)
that the generators $a$, $b$, $c$, $d$ admit a realization via
$q$-oscillators. Evidently, such a representation is not unique.

It is natural to realize $a$, $b$, $c$, $d_0$ as operators acting on
the $q$-analogue of the space $D(z_1,z_2)$. We shall denote this space
as $D_q(z_1,z_2)$.  The space $D_q(z_1,z_2)$ is spanned on the basic
vectors of form (remember that $[x]$ stands for $q$-numbers)
\beq{4.01} |j,m\rangle
=\frac{z_1^{j+m}z_2^{j-m}}{\sqrt{[j+m]![j-m]!}},\;\;\;\;j=
\ar{lllll}0,&\!\!\frac{1}{2},&\!\!1,&\!\!\frac{3}{2},&\!\! ...\er
\;\;\;\; m=-j,..,j.  \eeq One can define on $D_q(z_1,z_2)$ such a
scalar product that the system (\ref{4.01}) is orthonormal, that is $
<j,m|j^{\prime},m^{\prime}>=\delta_{j j^{\prime}} \delta_{m
m^{\prime}}$.

\rem This scalar product is a deformation of (\ref{3.3}). Its explicit
form makes use of the $q$-exponent and the Jackson integral. See
Ref.20 for details.

In all formulae concerning the space $D_q(z_1,z_2)$ we suppose that
$q$ is chosen as described in Sec.$\!$ II (i.e., it belongs either
to the real axis or to the unit circle at the complex plane).
In this case an
analogue of the rule of conjugation (\ref{3.5}) is
\beq{4.03}
(z_i)^{*}=z_i^{-1}\;[z_i\partial_i],\;\;\;\;\;
(z_i\partial_i)^{*}=z_i\partial_i.
\eeq

The formulae (\ref{4.105}) and (\ref{5.00}) imply that the generator
$b$ is a power of the operator of spin. Hence, on the space
$D_q(z_1,z_2)$ it is given by \beq{5.16} b = q^{z_1\,\partial_1 +
z_2\,\partial_2 + 1 } = q\,N_1\,N_2.  \eeq

Next, let us remind that we already know the limit versions of the
generators $a$, $b$, $c$, $d$ (see Proposition \ref{P6}; one should
take into account the rescaling (\ref{4.95})). Their appropriate
deformations for generic $q$ are described by
\begin{prop}\label{P15}
The set of operators (with arbitrary constants $\la_i$, $\nu_i$)
\beq{5.17} \ar{ll}
a=q^{\la_0}\,z_1^{-1/2}\,N_1^{\la_1}, & c=q^{\nu_0}\,z_1^{-1/2}\,z_2\,
N_1^{\la_1-2}\,N_2^{\nu_2}, \\ \\ b=q\,N_1\,N_2, &
d=-q^{\la_0-\nu_0+\nu_2}\,z_2^{-1}\,(N_2-N_2^{-1})\,N_1\,N_2^{-\nu_2} \er
\eeq 
satisfies (\ref{4.5})-(\ref{4.6}) and gives in the
limit $\ga\rar 0$ the generators found in (\ref{3.35b}).
\end{prop}

Although due to Theorem \ref{main} this proposition gives a family of 
representations for $\cal U$,
we again should impose an additional condition using the matrix
(\ref{L5}) as a standard (justification for this trick was given
above).
\begin{prop}\label{P16}
Matrix $L$ given by (\ref{4.9}) coincides with the matrix (\ref{L5})
taken in the representation (\ref{q33}) provided that \beq{5.19} b =
q\,N_1\,N_2,\;\;\;a\,c= q^{-1/2}\,z_1^{-1}\,z_2\,N_1^{-1/2}
\,N_2^{-1/2},\;\;\; a^{-1}\,c\,d_0=-N_1^{-1}\,N_2.  \eeq
\end{prop}

Comparing the statements of Propositions \ref{P15} and \ref{P16}, we
derive: 
\beq{5.20} 
\ar{ll} a=q^{\la_0}\,z_1^{-1/2}\,N_1^{3/4},\ &
b=q\,N_1\,N_2, \\
c=q^{\nu_0}\,
z_1^{-1/2}\,z_2\,N_1^{-5/4}\,N_2^{-1/2}, & \\
d_0=-q^{\la_0-\nu_0-1/2}\, z_2^{-1}\,N_1\,N_2^{3/2},\ &
q^{\la_0+\nu_0}= q^{-1/8} .  \er
\eeq
Substituting (\ref{5.20}) into (\ref{5.15}) (and remember that
(\ref{5.15}) is defined only up to a coefficient), we get \beq{5.21}
e^{i\xi} = q^{2\epsilon}\,z_1\,N_1^{\ga-1/2}\,N_2^{\ga-1}, \eeq where
$\ga$ and $\epsilon$ are arbitrary. Finally, substituting
(\ref{5.20})-(\ref{5.21}) into (\ref{5.9}), 
%% (one should be accurate when calculating $e^{\pm i\xi/2}$)
 we obtain (one should remember
 $U$ and $\wh{U}$ are defined only up to arbitrary scaling factor)
\beq{5.22} 
\wh{U} = \! \left( \ar{cc}
\frac{1}{\omega}\,\al_0\,z_1^{-1}\, N_1^{1-\ga/2}\,N_2^{3/2-\ga/2}
(N_1-N_1^{-1}) & \be_0 \,z_2\,N_1^{\ga/2-3/2} N_2^{\ga/2-1} \\ & \\ -
\frac{1}{\omega}\, q^{-1} \al_0\, z_2^{-1} N_1^{1/2-\ga/2}\,
N_2^{1-\ga/2} (N_2-N_2^{-1}) & \be_0\,z_1\,N_1^{\ga/2-1}
N_2^{\ga/2-1/2} \er \right).  
\eeq

It is easy to check that the family of matrices (\ref{5.22}) exactly
coincides with what was obtained in $q$-oscillator approach (see
Propositions \ref{P9} and  \ref{P10}).

\subsection*{\red \ns\bf D. Quantum Clebsch-Gordan coefficients. }
\red Using the connection formula (\ref{2.7}) we get from (\ref{5.22})
a family of matrices $U$ which provide possible representations of the
algebra $\cal U$. It is natural to study an action of the entries of
these matrices on the space $D_q(z_1,z_2)$ described above. On the
basic vectors (\ref{4.01}) these operators act as follows:
\beq{6.1}
 \ar{c} U_1\, |j,m\rangle = C_1\,q^{\frac{1}{2}(j-m+1)}
\sqrt{\frac{[j+m]}{[2j+1]}}\, |j-\frac{1}{2},m-\frac{1}{2}\rangle , \\
\\ U_2\, |j,m\rangle = C_2\,q^{-\frac{1}{2}(j+m)}
\sqrt{\frac{[j-m+1]}{[2j+1}]}\, |j+\frac{1}{2},m-\frac{1}{2}\rangle ,
\\ \\ U_3\, |j,m\rangle =-C_3\,q^{-\frac{1}{2}(j+m+1)}
\sqrt{\frac{[j-m]}{[2j+1]}}\, |j-\frac{1}{2},m+\frac{1}{2}\rangle , \\
\\ U_4\, |j,m\rangle =C_4\,q^{\frac{1}{2}(j-m)}
\sqrt{\frac{[j+m+1]}{[2j+1]}}\, |j+\frac{1}{2},m+\frac{1}{2}\rangle , \er
\eeq
where the coefficients $C_i$ do not depend on $m$.

Note that, similarly to the classical case, the operators $U_i$
correspond to the basic shifts on the model space. Comparing the
matrix elements $<j^{\prime},m^{\prime}|U_i|j,m>$ following from
(\ref{6.1}) with values of CGC for $U_q(sl(2))$ given by $q$-analogue
of the Van-der-Waerden formula,\REF{18,21} which for the
decomposition of $V_j\otimes V_{1/2}$ looks like following
\[ \left\{\ar{ccc} j & \frac 12 & j^{\prime \prime }  \\
m & m^{\prime } & m^{\prime \prime } \er \right\}_q =
\delta_{m^{\prime \prime},m + m^{\prime}}\; \left(\frac{[j+\frac
12-j^{\prime \prime }]![j+j^{\prime \prime }-\frac 12]![j^{\prime
\prime }+\frac 12-j]!}{[j+j^{\prime \prime }+\frac 32]!}
\right)^{1/2} \times \] \beq{6.2} \times q^{\frac{1}{2} (j+\frac{1}{2}
-j^{\prime \prime})(j+ j^{\prime \prime} +\frac{3}{2}) + jm^{\prime}
-\frac{1}{2} m} \times \eeq
\[ \times \sum\limits_{r\geq 0}\frac{(-1)^r q^{-r(j+j^{\prime \prime}+
\frac 32)} \left( [j+m]![j-m]![j^{\prime \prime}+m^{\prime\prime
}]![j^{\prime \prime }-m^{\prime \prime }]![2j^{\prime \prime
}+1]\right) ^{1/2} } {[r]![j+\frac 12-j^{\prime \prime
}-r]![j-m-r]![\frac 12+m^{\prime } -r]![j^{\prime\prime }-\frac
12+m+r]![j^{\prime \prime }-j-m^{\prime }+r]!},\] we establish the
following correspondence
\[ \! \ar{ll}
\langle j^{\prime \prime} , m^{\prime \prime} \mid U_1\mid j,
m\rangle= \delta_{j^{\prime \prime },j- \frac 12 }\,
\al_0\,q^{(1-\ga/2)j-1/2}\, \left\{\ar{ccc} j & \frac 12 & j^{\prime
\prime} \\ m & -\frac 12 & m^{\prime \prime} \er \right\}_q , \\
\langle j^{\prime \prime} , m^{\prime \prime} \mid U_2 \mid j,
m\rangle= \delta_{j^{\prime \prime },j+ \frac 12 }\,
\be_0\,q^{(\ga/2-1)j}\,\left\{\ar{ccc} j & \frac 12 & j^{\prime
\prime} \\ m & -\frac 12 & m^{\prime \prime} \er \right\}_q , \\
\langle j^{\prime \prime} , m^{\prime \prime} \mid U_3 \mid j,
m\rangle= \delta_{j^{\prime \prime },j- \frac 12 }\,
\al_0\,q^{(1-\ga/2)j-1/2}\, \left\{\ar{ccc} j & \frac 12 & j^{\prime
\prime} \\ m & \frac 12 & m^{\prime \prime} \er \right\}_q, \\ \langle
j^{\prime \prime} , m^{\prime \prime} \mid U_4 \mid j, m\rangle=
\delta_{j^{\prime \prime },j+ \frac 12 }\,
\be_0\,q^{(\ga/2-1)j}\,\left\{\ar{ccc} j & \frac 12 & j^{\prime
\prime} \\ m & \frac 12 & m^{\prime \prime}\er \right\}_q .  \er \]
Thus we derive an analogue of Proposition \ref{P7}:
\begin{prop}\label{P17}
The generators $U_i$ of the algebra $\cal U$ are operators of the basic 
shifts on the model space for $U_q(sl(2))$ and they generate the
$q$-Clebsch-Gordan coefficients corresponding to decomposition of the
product $V_j\otimes V_{\frac 12}$ of irreps of $U_q(sl(2))$.
\end{prop}

\rem Putting $\al_0=q^{1/2}$, $\be_0=1$ and $\ga=2$ in (\ref{5.22}),
we get the following generating matrix
\beq{6.3} 
U = \left( \ar{cc} z_1^{-1}\,
[z_1 \partial_1 ] q^{\frac 12 (z_2 \partial_2 +1)} & z_2 \, q^{-\frac
12 z_1 \partial_1 } \\ \\-z_2^{-1}\, [z_2 \partial_2 ] q^{-\frac 12
(z_1 \partial_1 +1)} & z_1 \, q^{\frac 12 z_2 \partial_2 } \er \right)
\,\frac{1}{\sqrt{\p}} \ , \;\;\;\;
{\rm Det} U = 
%%% (\wh{U}_1 \wh{U}_4 - \wh{U}_2 \wh{U}_3 )\frac{1}{\p} = 
q^{1/2} \ ,
\eeq 
which may be called ``exact'' as it satisfies (\ref{6.1}) with 
$C_i=1$. 
The question about unitarity of the matrix (\ref{6.3}) is discussed
in Appendix B.

\subsection*{\red \ns\bf E. Generalized Wigner-Eckart theorem.}

\red As we demonstrated in the previous section, entries of the matrix
$U_0$ are tensor operators of spin $1/2$ for ${\cal J}=sl(2)$, hence
they provide a realization of the Wigner-Eckart theorem. Let us now
consider the matrix $U$ from this point of view.

The theory of tensor operators for quantum algebras was discussed by
many authors (see, e.g., Ref.22). In particular, generalized Wigner-Eckart
theorem (in the case of $\J_q=U_q(sl(2))$ reads as follows.
\begin{theorem} \label{TWE}
Let $l_{+}$, $l_{-}$ and $l_{3}$ be the generators of $U_q(sl(2))$ and
let $T^{j\;}_{\;m}$, $m=-j,..,j$ be a set of operators acting on the
deformed model space $\cal M$ and obeying the commutation relations
\beq{WE1}
{} [l_3,T^{j\;}_{\;m}]=m\,T^{j\;}_{\;m},\;\;\;
{} l_{\pm}\,T^{j\;}_{\;m}\,q^{l_3} - q^{l_3\mp 1}\, T^{j\;}_{\;m}\,l_{\pm} 
=\sqrt{[j\mp m][j\pm m+1]}\,T^{j\;}_{\;m\pm 1}.
\eeq
Then the matrix elements of $T^{j\;}_{\;m}$ on $\cal M$ are proportional to
$q$-Clebsch-Gordan coefficients:
\[ <j^{\prime\prime} m^{\prime\prime} |T^{j\;}_{\;m}|j^{\prime} m^{\prime} >
= C_{j^{\prime\prime}j^{\prime}}^{j}\;\left\{\ar{ccc}j^{\prime} & j &
j^{\prime\prime}\\ m^{\prime} & m & m^{\prime\prime} \er \right\}_q,    \]
where the coefficients $C_{j^{\prime\prime}j^{\prime}}^j$ do not depend on
$m$, $m^{\prime}$, $m^{\prime\prime}$.
\end{theorem}

Proposition \ref{P17} implies that $U_i$ may be regarded as
q-tensor operators. Indeed, using (\ref{5.22}) and (\ref{q33}), one can 
check that $U_i$ satisfy (\ref{WE1}) (one obtains for $U_i$
deformations of relations (\ref{we20})). Similarly to the classical case 
we have the following.
\begin{prop}\label{P18}
The generators $U_i$ of the algebra ${\cal U}$ form a basis for $q$-tensor 
operators of spin $1/2$, 
that is components $T^{1/2}_{1/2}$ and 
$T^{1/2}_{-1/2}$ of any $q$-tensor operator of spin 1/2 can be realized as
linear combinations of $U_i$:
\beq{WE5}
T^{1/2}_{-1/2}= \mu(p) \;U_1+\nu(p) \;U_2,\;\;\; 
T^{1/2}_{1/2}=\mu(p) \;U_3+\nu(p) \;U_4,
\eeq
where $\mu(p)$ and $\nu(p)$ are functions only of $p$.
\end{prop}
\rem Unlike the classical case, solution (\ref{5.22}) gives a
family of matrices $U$. However, the corresponding matrix elements
$<j^{\prime\prime} m^{\prime\prime} |U_i|j^{\prime} m^{\prime} > $
differ only by factors which do not depend on $m^{\prime}$,
$m^{\prime\prime}$. Thus, any representative of obtained family of
matrices $U$ may be used in Proposition \ref{P18}.

Let us end the description of the algebra $\cal U$ from the point of view
of theory of $q$-tensor operators with the following statement:
\begin{prop}\label{P19}
The matrices $U$ and $L$ defined in the Theorem 1 obey the relation
\beq{WE10}
 R_- \,\up{1}{U}\,\up{2}{L} = \up{2}{L}\,R_+ \,\up{1}{U}.
\eeq
\end{prop}
The proof is straightforward
\[ \ar{c} R_- \,\up{1}{U}\,\up{2}{L} \, =
R_- \,\up{1}{U}\,\up{2}{U}\,\up{2}{D}\,(\up{2}{U})^{-1}= \,
\up{2}{U}\, \up{1}{U}\, {\cal R}_- \,\up{2}{D}\,(\up{2}{U})^{-1} = \\ [1mm]
=\, \up{2}{U}\, \up{1}{U}\,\up{2}{D}\,\sigma\, {\cal R}_+ \,(\up{2}{U})^{-1} 
= \, \up{2}{U}\,\up{2}{D}\,\up{1}{U}\, {\cal R}_+ \,(\up{2}{U})^{-1} =
\up{2}{U}\,\up{2}{D}\,(\up{2}{U})^{-1}\, R_+ \,\up{1}{U} = \,
\up{2}{L}\,R_+ \,\up{1}{U}; \er   \]
it makes use the relations (\ref{1.8})-(\ref{1.9}) and the property
(\ref{RD}).

A remarkable fact is that (\ref{WE10}) may be used for definition of
$q$-tensor operators instead of (\ref{WE1}). Indeed, in the limit 
$\ga\rar 0$ it turns into
\beq{WE3}
[\up{1}{U_0}, \up{2}{L}_0 ] = \Lambda\, \up{1}{U_0},\;\;\;\;
\Lambda=\left(\ar{cccc} 1/2& & & \\ &-1/2&1& \\ &1&-1/2& \\ & & &1/2
\er\right).
\eeq
Using the explicit form of $L_0$ given in (\ref{3.05}), one can easily
check that this matrix relation is equivalent to (\ref{we20}).
More on $R$-matrix description of $q$-tensor operators is given in Ref.23.

\subsection*{ \red \ns\bf  CONCLUSION.}

 In this paper we have constructed the $q$-analogue of the phase space
$T^*\B$ and clarified its role in description of the model representation
of the corresponding quantum group $G_q$. We unraveled a connection between
the algebras generating by entries of matrix $(A,B)$, $(U,D)$ and $(L,Q)$.
The general formulae were concretized by the example of $G=SL(2)$.

An extension of the described scheme to the case of arbitrary group $G$
will definitely improve understanding of the role played by the matrix
${\cal R}(p)$ which so far has been discussed in the literature much
less then the standard matrix $R$.

The results of this paper can be generalized in several directions even for
the case of $SL(2)$. The first is consideration of the matrix $U$ with an
auxiliary space corresponding to the higher spin representation. It must
lead to an exact form of generating matrix for all CGC. The work in this
direction is in progress now.
The second point to be discussed is the case of $q$ being a root of unity.
The structure of ${\cal R}(p)$ allows to hope that reduction on so-called
"good" representations will be quite natural in our formalism.
However, this case is to be examined more carefully.

\subsection*{\ns\bf \red Acknowledgments. }

We are grateful to A.Yu. Alekseev, P.P. Kulish and V. Schomerus for
stimulating discussions and useful comments. We would like to thank
Prof. A. Niemi for hospitality at TFT, University
of Helsinki, where this work was begun.

This work was partially supported by ISF grant R2H000 and by INTAS grant.

\subsection*{ \red \ns\bf  Appendix A: Proof of Theorem 1. }

Using (\ref{1.8})-(\ref{1.9}) together with the identity (\ref{RD})
and taking into account that matrices $\up{1}{D}$, $\up{2}{D}$ and
$\sigma$ mutually commute, we check
\[ \ar{c}
\up{1}{L} R_-^{-1} \up{2}{L} R_- = \,\up{1}{U} \up{1}{D} (\up{1}{U})^{-1}
R_-^{-1}\up{2}{U}\, \up{2}{D} (\up{2}{U})^{-1} R_- =\\  \\
=\,\up{1}{U} \up{1}{D}\,\up{2}{U} {\cal R}_-^{-1}(p)\, (\up{1}{U})^{-1}\! 
\up{2}{D} (\up{2}{U})^{-1} R_- = \, \up{1}{U}\, \up{2}{U} \up{1}{D} \sigma\,
{\cal R}_-^{-1}(p) \up{2}{D} \sigma\, 
(\up{1}{U})^{-1} (\up{2}{U})^{-1} R_- =\\  \\
=R_+^{-1} \up{2}{U}\, \up{1}{U} {\cal R}_+ (p) \up{1}{D}\sigma\,
{\cal R}_-^{-1}(p) \up{2}{D} \sigma\, {\cal R}_- (p)\, (\up{2}{U})^{-1} 
(\up{1}{U})^{-1} =\\  \\
=R_+^{-1}\up{2}{U}\,\up{1}{U} \up{2}{D} \sigma\up{1}{D} {\cal R}_- (p) \,
(\up{2}{U})^{-1} (\up{1}{U})^{-1} = R_+^{-1} \up{2}{U} \up{2}{D}\, \up{1}{U}
{\cal R}_+ (p)\, \sigma\up{1}{D}(\up{2}{U})^{-1} (\up{1}{U})^{-1} =\\  \\
=R_+^{-1} \up{2}{U}\up{2}{D}\, \up{1}{U} {\cal R}_+ (p)\, (\up{2}{U})^{-1}\!
\up{1}{D}(\up{1}{U})^{-1} = R_+^{-1} \up{2}{U} \up{2}{D} (\up{2}{U})^{-1}
R_+\up{1}{U} \up{1}{D} (\up{1}{U})^{-1} = R_+^{-1} \up{2}{L} R_+ \up{1}{L}.
\er\]

\subsection*{\red \ns\bf  Appendix B: 
 On conjugation of $\bf U_0$ and $\bf U$.  }

 {}First we consider the matrix $U_0$. Using the rules of
conjugation (\ref{3.5}) (and taking into account that $p^*=p$), one can
check that the matrix conjugated to $U_0$ does not coincide with
$U_0^{-1}$; that is, the matrix $U_0$ itself is not unitary. However,
it turns out
that the transposed matrix (one should remember that in general
$(U^T)^{-1}\neq (U^{-1})^T$ for matrices with non-commuting entries)
\[U_0^T = \left( \ar{cc} \partial_1 & -\partial_2\\
z_2 & z_1 \er\right)\,\sqrt{\frac{\hbar}{p}} \]
satisfies the unitarity condition:
\[ (U_0^T)^* = \sqrt{\frac{\hbar}{p}}\,\left(\ar{cc} z_1 & \partial_2 \\
-z_2 & \partial_1 \er \right) = (U_0^T)^{-1}. \]

In the deformed case (recall that $q$ can be either real 
or $|q|=1$) the matrix $U$ includes the operator $N$ which conjugates in 
different ways for the different choices of $q$. 
Let us consider the matrix $U$ given by (\ref{6.3}). The conjugated matrix 
can be constructed according to the rules (\ref{4.03}). Using the formula 
(\ref{q1}), one can check that the unitarity condition
$(U^T)^{*}\,U^T=\,U^T\,(U^T)^{*}= I $
(i.e., the same as in the non-deformed case) for the transposed matrix
holds only for real $q$. For $|q|=1$ see Ref.13.

\subsection*{\red \ns\bf Appendix C: Proof of Theorem 4. }

Here we complete the proof of Theorem \ref{main}, i.e., we have to
prove that matrix (\ref{5.9}) satisfies (\ref{2.8})-(\ref{2.11}) if
the following  conditions

\hspace*{2cm}$\ar{c}  \\  \al + \be = 0,\;\;
\ga+\be-\al+1=0,\;\; v(b)=1,\;\;w(b)=1/\omega \\  \\ \er
$\hfill (C\,1)

$\hspace*{-6mm}$are fulfilled. 

First, using relations (\ref{4.5}), (\ref{4.8}), (\ref{5.10}) and
conditions (C\,1), we check
\[\ar{l} \wh{U}_1 \wh{U}_2 = a\,(b+a^{-1}c\,d_0)\, e^{-i\xi}\;c\,e^{i\xi}=
q^{-1/2+\be/2} c\,a(b+ q\,a^{-1}c\,d_0) \,e^{i\xi} \, e^{-i\xi} =\\ [1.5mm]
= q^{\al/2+\be/2} c\,e^{i\xi}\,a(b+ q^{(\ga+\be-\al+1)/2}\,a^{-1}c\,d_0) \,
e^{-i\xi} =  \wh{U}_2 \wh{U}_1  . \er  \]
\[\ar{l} \wh{U}_1 \wh{U}_3 = a\,(b+a^{-1}c\,d_0)\,e^{-i\xi}\;c^{-1}\,
(b^{-1}+q^{-1}\,a^{-1}c\,d_0)\, e^{-i\xi} = \\ [1.5mm]
= q^{1/2-\be/2}c^{-1}\,a\,
(b+q^{-1}\,a^{-1}c\,d_0)\,e^{-i\xi}\,(b^{-1}+q^{-1}\,a^{-1}c\,d_0)\,
e^{-i\xi} =\\ [1.5mm] 
=q^{-\be/2}c^{-1}\,a\,(b^{-1}+q^{-1}\,a^{-1}c\,d_0)\,
(b+q^{-1}\,a^{-1}c\,d_0)\,e^{-i\xi}\,e^{-i\xi} = \\ [1.5mm]
= q^{-1/2-\be/2}c^{-1} \,(b^{-1}+q^{-1}\,a^{-1}c\,d_0)\,a\,
(b+q^{-1}\,a^{-1}c\,d_0)\,e^{-i\xi}\, e^{-i\xi} = \\ [1.5mm]
= q^{-1-\al/2-\be/2}c^{-1}\,(b^{-1}+q^{-1}\,a^{-1}c\,d_0)\,a\,
e^{-i\xi}\,(b+a^{-1}c\,d_0)\,e^{-i\xi} = q^{-1} \wh{U}_3 \wh{U}_1 .\er \]
The rest of relations (\ref{2.8}) can be proved similarly.

Next, note that relation (\ref{2.10}) can be rewritten as follows

\hspace*{3cm}$\ar{c} \\  \wh{U}_1\;\wh{U}_4\,(b-b^{-1})-\wh{U}_4\;\wh{U}_1\,
(q\,b-q^{-1} \,b^{-1}) =-\omega\, \wh{U}_3 \; \wh{U}_2 \, b. \\  \\ 
\er $ \hfill (C\,2)

$\hspace*{-7.5mm}$ To prove this equality we transform its l.h.s.
and r.h.s. as follows
\[\ar{l}\wh{U}_1\;\wh{U}_4\,(b-b^{-1})-\wh{U}_4\;\wh{U}_1\,(q\,b-q^{-1}
\,b^{-1}) =
 a\,(b+a^{-1}c\,d_0)\,e^{-i\xi} \; a^{-1} e^{i\xi}\,(b-b^{-1})- \\ [1.5mm]
- a^{-1} e^{i\xi} \; a\,(b+a^{-1}c\,d_0)\,e^{-i\xi} (q\,b-q^{-1}\,b^{-1}) =
 q^{-\al/2}\,(q^{1/2}\,b+q^{-1/2}\,a^{-1}c\,d_0)\,(b-b^{-1}) - \\  [1.5mm]
- q^{-\al/2}\,(q^{-1/2}\,b+q^{1/2}\,a^{-1}c\,d_0)\,(q\,b-q^{-1}\,b^{-1}) =
-q^{-\al/2}\,\omega\,(q^{1/2}\,a^{-1}c\,d_0\,b+q^{-1/2}); \\ \\
\wh{U}_3\;\wh{U}_2 \,b =c^{-1}\,(b^{-1}+q^{-1}\,a^{-1}c\,d_0)\,
e^{-i\xi}\,c\,e^{i\xi} \,b = \\  [1.5mm]
= q^{\be/2}\,(q^{-1/2}\,b^{-1}+q^{1/2}\,
a^{-1}c\,d_0)\,b = q^{\be/2}\,(q^{-1/2}+q^{1/2}\,a^{-1}c\,d_0\,b). \er  \]
Thus, the equality (C\,2) is fulfilled if conditions (C\,1) are valid.
The relation (\ref{2.11}) can be proved in the same way.

\subsection*{\red \ns\bf  References.}

%\footnotesize
\small

\REF{1}{I.N.Bernstein, I.M.Gelfand, S.I.Gelfand, \ {\it Models for
representations of Lie groups}. In: Proceedings of I.G.Petrovsky seminars,
{\bf 2}, 3 (1976) (in Russian);\
I.N.Bernstein, I.M.Gelfand, S.I.Gelfand, \ Funct. analysis and its
applications, {\bf 9}, No 4, 61 (1975).  } \\ [1mm]
\REF{2} {A.M.Polyakov, \ Mod.Phys.Lett. {\bf A 2}, 893 (1987).} \\[1mm]
\REF{3} {A.Yu.Alekseev, S.L.Shatashvili, \ Comm. Math. Phys. {\bf 128},
197 (1990); \\
\ H.La, P.Nelson, A.S.Schwarz, \
Comm. Math. Phys. {\bf 134}, 539 (1990).} \\ [1mm]
\REF{4} {L.D.Faddeev, \
Comm. Math. Phys. {\bf 132}, 131 (1990). } \\ [1mm]
\REF{5} {A.Yu.Alekseev, L.D.Faddeev, \ Comm. Math. Phys. {\bf 141}, 413
(1991).} \\ [1mm]
\REF{6} {F.Falceto, K.Gawedzki, \ J. Geom. Phys.
 {\bf 11}, 251 (1993).} \\ [1mm]
\REF{7} {A.A.Kirillov, \ {\it Elements of the theory of representations}.
Springer-Verlag (1979);\\
 \ A.Yu.Alekseev, L.D.Faddeev, S.L.Shatashvili, \
J. Geom. Phys. {\bf 5}, 391 (1989). } \\ [1mm]
\REF{8} {A.Connes,\ {\it Noncommutative geometry}.
  (Academic Press, 1994).} \\ [1mm]
\REF{9} {L.D.Faddeev, N.Yu.Reshetikhin, L.A.Takhtajan, \
Algebra and analysis {\bf 1}, 193 (1990). } \\ [1mm]
\REF{10} {A.Yu.Alekseev, L.D.Faddeev, M.A.Semenov-Tian-Shansky,
A.Yu.Volkov, \ {\it The unraveling of the quantum group structure in the
WZNW theory}, preprint CERN-TH-5981/91 (1991). } \\  [1mm]
\REF{11} {N.Yu.Reshetikhin, M.A.Semenov-Tian-Shansky, \ Lett. Math. Phys.
{\bf 19}, 133 (1990). } \\ [1mm]
\REF{12} {P.P.Kulish, N.Yu.Reshetikhin, \ Zapiski Nauch.Semin. LOMI
{\bf 101}, 101 (1981) (in Russian). } \\ [1mm]
\REF{13} {A.G.Bytsko, V.Schomerus,\ {\it Vertex operators -- from toy
model to lattice algebras}, preprint q-alg/9611010 (1996).} \\ [1mm]
\REF{14} {A.P. Isaev,\ J. Phys. {\bf A 29}, 6903 (1996).} \\  [1mm]
\REF{15} {J-L.Gervais, A.Neveu, \ Nucl. Phys.
  {\bf B 238}, 125 (1984);\\
\ E.Cremmer, J-L.Gervais, \ Comm. Math. Phys.
 {\bf 134}, 619 (1990). } \\ [1mm]
\REF{16} {O.Babelon, \ Comm. Math. Phys.
       {\bf 139}, 619 (1991).} \\ [1mm]
\REF{17} {J.Avan, O.Babelon, E.Billey,\ Comm. Math. Phys. {\bf 178},
281 (1996). } \\ [1mm]
\REF{18} {A.N.Kirillov, N.Yu.Reshetikhin, \ Adv. Series in Math. Phys.
{\bf 11}, 202 (World Scientific, 1990).}  \\ [1mm]
\REF{19} {A.O.Barut, R.Raczka, \ {\it Theory of group representations
and applications} (Scient. Publishers, 1977);\\  \
L.C.Biedenharn, J.D.Louck, \ {\it Angular momentum in quantum physics},
Encyclopedia of mathematics and its applications, {\bf v.8}
(Addison-Wesley, 1981).  \\ [1mm]
\REF{20} {M.Arik, D.D.Coon, \ J. Math. Phys. {\bf 17}, 524 (1976);\\
\ L.C.Biedenharn, \ J. Phys. {\bf A 22}, L873 (1989);\\ \
A.J.Macfarlane, \ J. Phys. {\bf A 22}, 4581 (1989);\\  \
P.P.Kulish, E.Damaskinsky, \ J. Phys. {\bf A 23}, L415 (1990);\\ \
M.Chaichian, P.P.Kulish, \ Phys. Lett. {\bf B 234}, 72 (1990).} \\ [1mm]
\REF{21} {M.Nomura. \ J. Math. Phys. {\bf 30}, 2397 (1990);\\ \
L.Vaksman. \ Sov. Math. Dokl. {\bf 39}, 467 (1989);\\ \
H. Ruegg. \ J. Math. Phys. {\bf 31}, 1085 (1990).  } \\ [1mm]
\REF{22} {L.C.Biedenharn, M.Tarlini, \ Lett. Math. Phys. {\bf 20},
271 (1990);\\  \
K.Bragiel, \, Lett. Math. Phys. {\bf 21}, 181 (1991);\\  \
G.Mack, V.Schomerus,\, Phys. Lett. {\bf B 267}, 207 (1991);\\ \
V.Rittenberg, M.Scheunert, \, J. Math. Phys.
 {\bf 33}, 436 (1992).} \\ [1mm]
\REF{23} {A.G.Bytsko,\ {\it Tensor operators in $R$-matrix approach},
preprint DESY-95-254 (1995).}

\end{document}